# Robust Spectral Analysis

Andreas Hagemann*

*University of Notre Dame*

August 26, 2013

In this paper I introduce quantile spectral densities that summarize the cyclical behavior of time series across their whole distribution by analyzing periodicities in quantile crossings. This approach can capture systematic changes in the impact of cycles on the distribution of a time series and allows robust spectral estimation and inference in situations where the dependence structure is not accurately captured by the auto-covariance function. I study the statistical properties of quantile spectral estimators in a large class of nonlinear time series models and discuss inference both at fixed and across all frequencies. Monte Carlo experiments and an empirical example illustrate the advantages of quantile spectral analysis over classical methods when standard assumptions are violated.

## 1. Introduction

Classical spectral analysis uses estimates of the spectrum or spectral density, a weighted sum of auto-covariances, to quantify the relative magnitude and frequency of cycles present in a time series. However, if the dependence structure is not accurately captured by the auto-covariance function, for example, because the time series under consideration is uncorrelated or heavy-tailed, then spectral analysis can provide only uninformative or even misleading results. In this paper I discuss estimation and inference for a new class of spectral densities that summarize the cyclical behavior across the whole distribution of a time series by analyzing how frequently a process crosses its marginal quantiles. Functions from this class, which I refer to as quantile spectra or quantile spectral densities, are similar to classical spectral densities in both shape and interpretation, but can capture systematic changes in the impact of cycles on the distribution of a time series. Such changes arise naturally in a variety of modern time series models, including stochastic volatility and random coefficient autoregressive models, and cannot be identified through classical spectral analysis, where cycles are assumed to be global phenomena with a constant effect on the whole distribution. Quantile spectral analysis fundamentally changes this view because it distinguishes between the effects of cycles at different points of the distribution of a process and permits a local focus on the parts of the distribution that are most affected by the cyclical structure.

---

*Department of Economics, University of Notre Dame, 434 Flanner Hall, Notre Dame, Indiana 46556, USA. Tel.: +1 (574) 631-1688. Fax: +1 (574) 631-4783. E-mail address: ahageman@nd.edu. I would like to thank Roger Koenker for his continued support and encouragement. I would also like to thank Dan Bernhardt, Xiaofeng Shao, and Rolf Tschernig for comments and discussions. I gratefully acknowledge financial support from a UIUC Summer Research Fellowship. All errors are my own.

Spectral analysis has traditionally played an important role in the analysis of economic time series; see, among many others, Granger (1966), Sargent (1987, chap. 11), Diebold, Ohanian, and Berkowitz (1998), and Qu and Tkachenko (2011), where the shape of the sample spectral density is typically taken to be one of the "stylized facts" that the predictions of a model must match. For macroeconomic data, these stylized facts often refer to high-frequency (seasonal) and low-frequency (business cycle) peaks in the spectrum. However, both observed data and the posterior distributions of economic models can exhibit heavy tails (Cogley and Sargent, 2002) that can induce peaks at random in the sample spectra of the data and the model output, invalidating comparisons between the two. For financial data, the stylized facts include the absence of auto-correlation, i.e., peakless spectra, and heavy-tailed marginal distributions (Cont, 2001). Stochastic volatility models such as GARCH processes (Bollerslev, 1986) can cross almost every quantile of their distribution in a periodic manner and at the same time satisfy these and other stylized facts, leading the researcher to incorrectly conclude from the spectrum that no periodicity is present. Bispectra and higher-order spectra can possibly detect cycles in quantile crossings, but rely on the presumption of light tails since they require the existence of at least third moments to be well defined and sixth moments to be estimated reliably (see, e.g., Rosenblatt and Van Ness, 1965). Financial time series such as log-returns of foreign exchange rates or stock prices may lack finite fourth or even third moments (Loretan and Phillips, 1994; Longin, 1996).

My proposed approach is robust to each of these concerns: Quantile spectral methods consistently recover the spectral shape and detect periodicities even in uncorrelated or heavy-tailed data. Inference about quantile spectra both at fixed frequencies and across frequencies does not require assumptions about the moments of the process. Although moment conditions can be used to verify some of the assumptions below, arbitrarily low fractional moments suffice. Because several common time series models can induce situations where cycles are present at some but not at all quantiles, I also provide a general Cramér-von Mises specification test for peakless quantile spectra. Under conditions that are routinely imposed in the literature when testing for the absence of peaks, these tests are distribution-free and, depending on the strength of the assumptions, sometimes even exact in finite samples. The test remains valid asymptotically under much weaker conditions when a bootstrap approximation is used.

Several recent papers apply quantiles in spectral or correlogram (auto-correlation) analy-sis. Li (2008, 2012) obtains robust spectral estimators via quantile regressions for harmonic regression models. Although his estimation method is quite different from that developed here, there is some overlap in our results. I provide a detailed discussion in section 3. Katkovnik (1998) relies on the same idea as Li (2008), but only works with sinusoidal models and iid noise. Linton and Whang (2007) introduce the "quantilogram," a correlogram that is essentially the inverse Fourier transform of a quantile spectrum, but their focus is on testing for directional predictive ability of financial data in the time domain, rather than spectral analysis. Chung and Hong (2007) test for directional predictive ability with the generalized spectrum (Hong, 1999) by investigating the frequency domain behavior of processes around a given threshold. This approach is similar in spirit to quantile spectral analysis but, as Linton and Whang point out, Chung and Hong rescale their data with sample standard deviations but do not account for the randomness introduced by the rescaling in the derivation of their tests. In contrast, the scaling of the data for quantile spectral analysis is provided automatically through the marginal quantile function and all of my results are derived under the assumption that these quantiles are estimated.

Other robust spectral methods are discussed by Kleiner and Martin (1979) and Klüppelburg and Mikosch (1994): Kleiner and Martin focus on time series where the dependence structure is accurately captured by an autoregressive model of sufficiently high order. Quantile spectral analysis differs from these methods in that it is completely nonparametric and, most importantly, it robustly estimates cyclical components even when an autoregression is not an appropriate model for the data. Klüppelburg and Mikosch robustify classical spectral estimates by a self-normalization procedure to estimate normalized spectra under arbitrarily weak fractional moments conditions. However, their results are of limited use for applications because little is known about the asymptotic distribution of their procedures. In contrast, I show that quantile spectral estimates have relatively simple asymptotic distributions even when no moments exist.

After completing the first draft of this manuscript (Hagemann, 2011), the papers by Dette, Hallin, Kley, and Volgushev (2011) and Lee and Subba Rao (2011) became available. Both describe methods based on analyzing cross-covariances of quantile hits and copulas in the frequency domain that are similar to the quantile spectral estimators presented here. However, both Dette et al. and Lee and Subba Rao develop their methods as alternatives to the generalized spectrum to discover the presence of *any* type of dependence structure in time series data. My estimators are constructed to identify cyclical dependence.

The remainder of the paper is organized as follows: Section 2 discusses quantile spectral analysis and introduces two classes of estimators. Section 3 establishes the asymptotic validity of the estimators under weak regularity conditions. In section 4, I show the consistency of Cramér-von Mises tests for peakless spectra. The Monte Carlo experiments and an empirical example in section 5 illustrate the finite sample properties of the estimators and tests. Section 6 concludes. The Appendix contains auxiliary results and proofs.

I use the following notation throughout the paper: $1\{\cdot\}$ and $1_{\{\cdot\}}$ both represent the indicator function and $\|X\|_p$ denotes $(\mathbb{E}|X|^p)^{1/p}$; $\|\cdot\|$ abbreviates $\|\cdot\|_2$. Limits are as $n \to \infty$ unless otherwise noted and convergence in distribution is indicated by $\rightsquigarrow$. The inner product $\langle \cdot, \cdot \rangle_\Pi$ and norm $\|\cdot\|_\Pi$ are defined at the beginning of section 4.

## 2. Quantile Spectra and Two Estimators

This section introduces quantile spectral density estimation as a robust complement to classical spectral methods. Spectral analysis aims to reveal periodic behavior in a stationary time series $X_t$ with auto-covariance function $\gamma_X(j) := \mathbb{E}X_0 X_j - (\mathbb{E}X_0)^2$ at lag $j$ by estimating the *spectrum* or *spectral density* at frequency $\lambda$, defined as

$$f_X(\lambda) = \frac{1}{2\pi} \sum_{j \in \mathbb{Z}} \gamma_X(j) \cos(j\lambda), \qquad \lambda \in (-\pi, \pi]. \tag{2.1}$$

The auto-covariance function is typically taken to be absolutely summable to ensure that $f_X$ is continuous and symmetric about 0; a stochastic process that does not possess at least finite second moments cannot be meaningfully analyzed by the spectrum. If $f_X$ has a peak at $\lambda$, then $X_t$ is expected to repeat itself on average after $2\pi/\lambda$ units of time; for example, a monthly time series with a peak in the spectrum at $2\pi/3$ has a three-month cycle, with a higher value of $f_X$ corresponding to a more pronounced cycle. The primary goal of this paper is to develop spectral methods that go a step beyond summarizing the average impact of cycles by distinguishing between the effects of cycles at different points of the distribution of $X_t$.

The central idea is that if a stationary process $(X_t)_{t \in \mathbb{Z}}$ contains cycles, then its realizations will tend to stay above or below a given threshold in an approximately periodic manner. The pattern in which the process crosses a threshold at the center of its distribution reflects the most prominent cycles, but provides little information about their relative sizes. Patterns in threshold crossings near the extremes of the distribution help to identify amplitudes of these cycles and also recover periodicities that are obscured at the center of the distribution. The quantiles of $X_t$, arising from the quantile function $\xi_0(\tau) := \inf\{x : \mathbb{P}(X_0 \leq x) \geq \tau\}$, are natural choices for such thresholds because they give precise meaning to the notion of the center and extremes of a distribution. Spectral analysis of quantile crossing patterns can then discover cycles in the process and reveal the extent to which they are present at a given quantile without relying on moments.

To formalize this idea, pick probabilities $\tau \in (0, 1)$ corresponding to the marginal quantiles $\xi_0(\tau)$ of $X_t$. The variable of interest for the analysis is

$$V_t(\tau, \xi) = \tau - 1\{X_t < \xi\}, \qquad (\tau, \xi) \in (0, 1) \times \mathbb{R},$$

such that $V_t(\tau) := V_t(\tau, \xi_0(\tau))$ takes on the value $\tau - 1$ if $X_t$ is below its $\tau$-th quantile at $t$, and $\tau$ otherwise. Here the quantiles are not assumed to be known, which enables the researcher to choose the values of $\tau$ according to the amount and nature of information that is needed about the cyclical structure of the time series. For example, $\tau = 0.5$ only analyzes fluctuations about the median, whereas varying $\tau$ between 0.5 and 0.9 also provides information about the positive amplitudes by including values in the upper tail of the process.

If the distribution function of $X_t$ is continuous and increasing at $\xi_0(\tau)$, then the $\tau$-th quantile crossing indicator $V_t(\tau)$ is a bounded, stationary, mean-zero random variable with auto-covariance function $r_\tau(j) := \gamma_{V(\tau)}(j) = \mathbb{E} V_0(\tau) V_j(\tau)$. Periodicities in $V_t(\tau)$ are summarized by peaks in its spectral density

$$g_\tau(\lambda) := f_{V(\tau)}(\lambda) = \frac{1}{2\pi} \sum_{j \in \mathbb{Z}} r_\tau(j) \cos(j\lambda), \qquad (2.2)$$

which I refer to as the *$\tau$-th quantile spectrum* or *$\tau$-th quantile spectral density* in the sequel. Analyzing $g_\tau$ across a grid of probabilities $\tau \in (0,1)$ therefore reveals cycles in events of the form $\{X_t < \xi_0(\tau) : t \in \mathbb{Z}\}$, which in turn summarize $(X_t)_{t \in \mathbb{Z}}$ with arbitrary precision as long as the grid is fine enough.

As the next two examples show, quantile spectral analysis can in fact yield additional insights beyond classical spectral analysis; Linton and Whang (2007) consider similar models. I discuss estimation of quantile spectra below.

**Example 2.1** (Stochastic volatility)**.** Let $(\varepsilon_t)_{t \in \mathbb{Z}}$ be a sequence of iid mean-zero random variables and suppose the data are generated by the stochastic volatility model $X_t = \xi_0(\tau_0) + \varepsilon_t v(\varepsilon_{t-1}, \varepsilon_{t-2}, \dots)$ for some $\tau_0 \in (0,1)$, where $v > 0$ is a measurable function that drives the volatility of the process. If $X_t$ has finite second moments, then it is an uncorrelated time series and its spectrum contains no information about the dependence structure beyond that it is "flat," i.e., $f_X(\lambda) = \gamma_X(0)/(2\pi)$ at all frequencies.

Similarly, any stationary time series with a continuous and increasing distribution function at $\xi_0(\tau)$ satisfies $r_\tau(0) = \tau(1-\tau)$ and the stochastic volatility process also has the property that

$$r_{\tau_0}(j) = \mathbb{E} V_0(\tau_0) \big( \tau_0 - \mathbb{P}(X_j < \xi_0(\tau_0) \mid \varepsilon_{j-1}, \dots) \big) = \big( \tau_0 - \mathbb{P}(\varepsilon_j < 0) \big) \mathbb{E} V_0(\tau_0) = 0$$

for all $j > 0$. Therefore its $\tau_0$-th quantile spectrum will also flat in the sense that $g_{\tau_0}(\lambda) \equiv \tau_0(1 - \tau_0)/(2\pi)$. However, the other quantile spectra of the stochastic process will be informative because $r_\tau(j)$ generally does not vanish for $\tau \neq \tau_0$. □

**Example 2.2** (QAR)**.** Now let $(\varepsilon_t)_{t \in \mathbb{Z}}$ be a sequence of independent Uniform(0, 1) variables and consider the second-order quantile autoregressive (QAR(2)) process of Koenker and Xiao (2006),

$$X_t = \beta_0(\varepsilon_t) + \beta_1(\varepsilon_t) X_{t-1} + \beta_2(\varepsilon_t) X_{t-2} = \mathbb{E}(\beta_1(\varepsilon_0)) X_{t-1} + \mathbb{E}(\beta_2(\varepsilon_0)) X_{t-2} + Y_t,$$

where $Y_t = \beta_0(\varepsilon_t) + [\beta_1(\varepsilon_t) - \mathbb{E}(\beta_1(\varepsilon_0))] X_{t-1} + [\beta_2(\varepsilon_t) - \mathbb{E}(\beta_2(\varepsilon_0))] X_{t-2}$. Here $\beta_0$, $\beta_1$, and $\beta_2$ are unknown functions that satisfy regularity conditions that ensure stationarity and $X_t$ is assumed to be increasing in $\varepsilon_t$ conditional on $X_{t-1}, X_{t-2}$. Provided that its second moments exist, the sequence $(Y_t)_{t \in \mathbb{Z}}$ has no influence on the shape of the spectrum because it is an uncorrelated sequence that is also uncorrelated with the other variables on the right-hand side of the preceding display (Knight, 2006). Hence, if $\mathbb{E}\beta_1(\varepsilon_0) = \mathbb{E}\beta_2(\varepsilon_0) = 0$, the spectrum of $X_t$ satisfies $f_X(\lambda) \equiv \gamma_Y(0)/(2\pi)$ and classical spectral analysis cannot reveal anything about cycles in $X_t$. If $\mathbb{E}\beta_1(\varepsilon_0)$ and $\mathbb{E}\beta_2(\varepsilon_0)$ are nonzero, then the spectrum of $X_t$ is the same as that of an AR(2) process with the same mean and covariances as the QAR(2). If, instead, there is some $\tau_0 \in (0,1)$ such that $\beta_1(\tau_0) = \beta_2(\tau_0) = 0$, then the $\tau_0$-th quantile spectrum is also flat (see Example 4.3 below), but cycles can be recovered at other quantiles. Further, the quantile spectra of the QAR(2) process and those of an AR(2) process with the same mean and covariance structure will generally be different. □

For a given sample $\mathcal{S}_n := \{X_t : t = 1, \ldots, n\}$, I consider two estimators of the quantile spectrum that correspond to the periodograms and smoothed periodograms used in classical spectral analysis. The key difference from the classical case is that the variable of interest $V_t(\tau)$ is indexed by the unknown quantity $\xi_0(\tau)$ and therefore itself has to be estimated. To this end, let $\hat{\xi}_n(\tau)$ be the $\tau$-th sample quantile determined implicitly by solutions to the minimization problem

$$\min_{x \in \mathbb{R}} \sum_{t=1}^{n} \rho_\tau(X_t - x),$$

where $\rho_\tau(x) := x(\tau - 1\{x < 0\})$ is the Koenker and Bassett (1978) check function. Let $\hat{V}_t(\tau) := V_t(\tau, \hat{\xi}_n(\tau))$ be the estimate of $V_t(\tau)$. The $\tau$-*th quantile periodogram* is then the "plug-in" estimator

$$Q_{n,\tau}(\lambda) := \frac{1}{2\pi}\left|n^{-1/2}\sum_{t=1}^{n}\hat{V}_t(\tau)e^{-it\lambda}\right|^2 = \frac{1}{2\pi}\sum_{|j|<n}\hat{r}_{n,\tau}(j)\cos(j\lambda), \qquad (2.3)$$

where $i := \sqrt{-1}$ and $\hat{r}_{n,\tau}(j) := n^{-1}\sum_{t=|j|+1}^{n}\hat{V}_t(\tau)\hat{V}_{t-|j|}(\tau)$ for $|j| < n$. As I will show in the next section, the quantile periodogram inherits the properties of the classical periodogram in the sense that it allows the construction of valid confidence intervals, but does not provide consistent estimates for the spectrum of interest.

Consistent estimation of the quantile spectrum requires additional smoothing to assign less weight to the imprecisely estimated auto-covariances with lags $|j|$ near $n$. For this I apply the Parzen (1957) class of kernel spectral density estimators to the present framework. The estimators, which I refer to as *smoothed $\tau$-th quantile periodograms*, are given by

$$\hat{g}_{n,\tau}(\lambda) = \frac{1}{2\pi}\sum_{|j|<n} w(j/B_n)\hat{r}_{n,\tau}(j)\cos(j\lambda), \qquad (2.4)$$

where $B_n$ is a scalar "bandwidth" parameter that grows with $n$ at a rate specified in Theorem 3.6 below, and $w$ is a real-valued smoothing weight function from the set

$$\mathcal{W} := \Big\{w \text{ is bounded and continuous}, w(x) = w(-x) \; \forall x \in \mathbb{R},$$
$$w(0) = 1, \bar{w}(x) := \sup_{y \geq x}|w(y)| \text{ satisfies } \int_0^\infty \bar{w}(x)\,dx < \infty,$$
$$W(\lambda) := (2\pi)^{-1}\int_{-\infty}^\infty w(x)e^{-ix\lambda}\,dx \text{ satisfies } \int_{-\infty}^\infty |W(\lambda)|\,d\lambda < \infty\Big\}.$$

In the literature, $w$ and $W$ are usually called the *lag window* and *spectral window*, respectively. Both functions are also often referred to as kernels, although $w$ does not necessarily integrate to one.

*Remarks.* 1. The class $\mathcal{W}$ includes most of the kernels that are used in practice, for example the Bartlett (i.e., triangular), Parzen, Tukey-Hanning, Daniell, and quadratic-spectral windows. However, it excludes the truncated (also known as rectangular or Dirichlet) window. See Andrews (1991) and Brockwell and Davis (1991, pp. 359-362) for thorough descriptions of these windows and their properties. I provide a brief discussion on how to choose $w$ and $B_n$ at the end of the next section.

2. The restriction $\int_0^\infty \bar{w}(x)\,dx < \infty$ is not standard in the spectral density estimation literature. As pointed out by Jansson (2002), it is needed for asymptotic bounds on expressions such as $B_n^{-1} \sum_{|j|<n} |w(j/B_n)|$ that typically arise in consistency proofs of spectral density estimates indexed by estimated parameters; see also Robinson (1991).

3. Spectra are non-negative. It is therefore common practice to choose a window such that $W \geq 0$ to ensure non-negativity of the spectral density estimate; see, e.g., Andrews (1991) and Smith (2005). The condition $\int_{-\infty}^\infty |W(\lambda)|\,d\lambda < \infty$ is immediately satisfied for such windows in view of the inverse Fourier transform $w(x) = \int_{-\infty}^\infty e^{ix\lambda} W(\lambda)\,d\lambda$, which implies $\int_{-\infty}^\infty W(\lambda)\,d\lambda = w(0) = 1$ for $w \in \mathcal{W}$. The Tukey-Hanning window is an example of a window that satisfies $\int_{-\infty}^\infty |W(\lambda)|\,d\lambda < \infty$, but not $W \geq 0$.

The next section characterizes the asymptotic properties of the quantile and smoothed quantile periodograms.

## 3. Asymptotic Properties of Quantile and Smoothed Quantile Periodograms

In this section I construct confidence intervals for the quantile spectrum and establish the consistency of the smoothed quantile periodogram under regularity conditions. I also compare the quantile periodogram to the periodograms of Li (2008, 2012).

Throughout the remainder of the paper I assume that $(X_t)_{t\in\mathbb{Z}}$ is a nonlinear process of the form

$$X_t = Y(\varepsilon_t, \varepsilon_{t-1}, \varepsilon_{t-2}, \dots), \tag{3.1}$$

where $(\varepsilon_t)_{t\in\mathbb{Z}}$ is a sequence of iid copies of a random variable $\varepsilon$ and $Y$ is a measurable, possibly unknown function that transforms the input $\mathcal{F}_t := (\varepsilon_t, \varepsilon_{t-1}, \dots)$ into the output $X_t$. The class (3.1) includes a large number of commonly-used stationary time series models. For instance, the processes in Examples 2.1 and 2.2 are of this form; I provide other examples below Proposition 3.1 in this section.

The essential conditions for the estimation of spectra are restrictions on the memory of the time series. As pointed out by Wu (2005), for time series of the form (3.1) such restrictions are most easily implemented by comparing $X_t$ to a slightly perturbed version of itself. Let $(\varepsilon_t^*)_{t\in\mathbb{Z}}$ be an iid copy of $(\varepsilon_t)_{t\in\mathbb{Z}}$, so that the difference between $X_t$ and $X'_t := Y(\varepsilon_t, \ldots, \varepsilon_1, \varepsilon_0^*, \varepsilon_{-1}^*, \ldots)$ are the inputs before time $t = 1$. Define $\mathcal{X}_\tau(\delta) := \{\xi \in \mathbb{R} : |\xi_0(\tau) - \xi| \leq \delta\}$ and assume the following:

**Assumption A.** *For a given $\tau \in (0,1)$, there exist $\delta > 0$ and $\sigma \in (0,1)$ such that*

$$\sup_{\xi \in \mathcal{X}_\tau(\delta)} \|1\{X_n < \xi\} - 1\{X'_n < \xi\}\| = O(\sigma^n).$$

Intuitively, this condition requires the probability that $X_n$ is below but $X'_n$ is above a given threshold (or vice versa) to be sufficiently small for large $n$ as long as the threshold is near $\xi_0(\tau)$. It is the only dependence condition needed to construct asymptotically valid confidence intervals for the quantile spectrum. Assumption A avoids restrictions on the summability of the cumulants (Brillinger, 1975, pp. 19-21) of $X_t$ that are routinely imposed in the spectral estimation literature; see Andrews (1991) and the references therein. Cumulant conditions or "mixing" assumptions (Rosenblatt, 1984) that imply such conditions are sometimes difficult to establish for a given time series model and can easily fail or put unwanted restrictions on the parameter space when $X_t$ is, for example, generated by a standard GARCH process (Bollerslev, 1986).

Assumption A does not require the existence of any moments of $X_t$, but can be verified for most commonly-used stationary time series models at the expense of an arbitrarily weak moment restriction via the geometric moment contracting (GMC) property introduced by Hsing and Wu (2004). A time series of the form (3.1) is said to be GMC for some $\alpha > 0$ if $\|X_n - X'_n\|_\alpha = O(\varrho^n)$ for some $\varrho \in (0,1)$, where $\varrho$ may depend on $\alpha$.

**Proposition 3.1.** *Assumption A is satisfied if $F_X(x) := \mathbb{P}(X_0 \leq x)$ is Lipschitz continuous in a neighborhood of $\xi_0(\tau)$ and $\|X_n - X'_n\|_\alpha = O(\varrho^n)$ for some $\alpha > 0$ and $\varrho \in (0,1)$.*

The GMC property is satisfied for stationary (causal) ARMA, ARCH (Engle, 1982), GARCH, ARMA-ARCH, ARMA-GARCH, asymmetric GARCH (Ding, Granger, and Engle, 1993; Ling and McAleer, 2002), generalized random coefficient autoregressive (Bougerol and Picard, 1992), and QAR models; see Shao and Wu (2007) and Shao (2011b) for proofs and more examples. By Proposition 3.1, all of these models are included in the analysis if $F_X$ is Lipschitz near $\xi_0(\tau)$—a condition that is also needed for all of my results.

In addition to Lipschitz continuity, a restriction on $F_X$ is required to ensure both that $V_t(\tau)$ can be estimated consistently and that $\sqrt{n}(\hat{\xi}_n(\tau) - \xi_0(\tau))$ is bounded in probability:

**Assumption B.** *$F_X$ is Lipschitz continuous in a neighborhood of $\xi_0(\tau)$ and has a positive and continuous (Lebesgue) density at $\xi_0(\tau)$.*

This assumption, or slight variations thereof, is standard in the quantile estimation and regression literature; see, e.g., Koenker (2005, p. 120) and Wu (2007).

As a preliminary step towards inference about quantile spectra, the following result establishes the joint asymptotic distribution of the quantile periodogram on a subset of the *natural frequencies* $\ldots, -4\pi/n, -2\pi/n, 0, 2\pi/n, 4\pi/n, \ldots \subset (-\pi, \pi]$. More precisely, Theorem 3.2 shows that the usual convergence of the periodogram at different frequencies to independent exponential variables is not affected by the presence of the estimated quantile $\hat{\xi}_n(\tau)$.

**Theorem 3.2.** *Suppose Assumptions A and B hold for some $\tau \in (0,1)$. Let $\lambda_n = 2\pi j_n/n$ with $j_n \in \mathbb{Z}$ be a sequence of natural frequencies such that $\lambda_n \to \lambda \in (0,\pi)$ with $g_\tau(\lambda) > 0$. Then, for any fixed $k \in \mathbb{Z}$, the collection of quantile periodograms*

$$Q_{n,\tau}(\lambda_n - 2\pi k/n), Q_{n,\tau}(\lambda_n - 2\pi(k-1)/n), \ldots, Q_{n,\tau}(\lambda_n + 2\pi k/n)$$

*converges jointly in distribution to independent exponential variables with mean $g_\tau(\lambda)$.*

*Remarks.* 1. The natural frequencies induce invariance to centering in the quantity inside the modulus in (2.3) so we can write

$$n^{-1/2} \sum_{t=1}^{n} \hat{V}_t(\tau) e^{-it\lambda_n} = -n^{-1/2} \sum_{t=1}^{n} \Big(1\{X_t < \hat{\xi}_n(\tau)\} - F_X\big(\hat{\xi}_n(\tau)\big)\Big) e^{-it\lambda_n}.$$

Given the invariance, the strategy for the proof is to show that the empirical process on the right-hand side of the preceding display is stochastically equicontinuous with respect to an appropriate semi-metric on bounded sets near $\xi_0(\tau)$. For this I extend Andrews and Pollard's (1994) functional limit theorems to time series of the form (3.1) that satisfy Assumption A. The equicontinuity property and a result of Shao and Wu (2007) on classical periodograms at natural frequencies then yield the desired results.

2. If a quantile of interest $\xi_0(\tau)$ is assumed to be known, for example $\xi_0(0.5) = 0$ as in Li (2008), then Theorem 3.2 remains valid when (i) $\xi_0(\tau)$ is used in $Q_{n,\tau}$ instead of $\hat{\xi}_n(\tau)$, (ii) Assumption B is replaced by the condition that $F_X$ is continuous and increasing at $\xi_0(\tau)$, and (iii) Assumption A is replaced by Assumption C below with $\delta = 0$. This is a direct consequence of Shao and Wu's (2007) Corollary 2.1.

Theorem 3.2 yields a convenient way to construct point-wise confidence intervals for the quantile spectrum. The proof follows immediately from the properties of independent exponential variables. Example 3.4 provides an application.

**Corollary 3.3.** *Suppose the conditions of Theorem 3.2 are satisfied. Define $\bar{Q}_{n,\tau}(\lambda, k) = \sum_{|j| \le k} Q_{n,\tau}(\lambda_n + 2\pi j/n)/(2k+1)$, and let $\chi^2_{4k+2,\alpha}$ be the $\alpha$-th quantile of a $\chi^2$ distribution with $4k+2$ degrees of freedom. Then, for every fixed $k \in \mathbb{Z}$, the probability of the event*

$$g_\tau(\lambda) \in \left( \frac{(4k+2)\bar{Q}_{n,\tau}(\lambda,k)}{\chi^2_{4k+2,1-\alpha/2}}, \frac{(4k+2)\bar{Q}_{n,\tau}(\lambda,k)}{\chi^2_{4k+2,\alpha/2}} \right)$$

*converges to $1-\alpha$.*

**Example 3.4** (Testing for periodicities)**.** The processes in Examples 2.1 and 2.2 are instances where $V_t(\tau_0)$ is a white noise series for some $\tau_0 \in (0,1)$. Then the $\tau_0$-th quantile spectrum of $X_t$ is $\tau_0(1-\tau_0)/(2\pi)$ at all frequencies and therefore contains no periodicities at that quantile. Because a spike in the periodogram could either be evidence for a periodicity or an artifact generated by the sample, this leads to the problem of testing whether the $\tau_0$-th quantile spectrum behaves like a flat quantile spectrum at a given frequency.

By Corollary 3.3, this hypothesis can be rejected at level $\alpha$ if the confidence interval in the Corollary does not contain $\tau_0(1-\tau_0)/(2\pi)$. The same type of test is not as simple in classical spectral analysis because (2.1) reduces to the unknown quantity $\gamma_X(0)/(2\pi)$ if $X_t$ is white noise. I extend the idea of testing for flatness in section 4 to provide a test for the more general hypothesis that $g_{\tau_0}(\lambda) = \tau_0(1-\tau_0)/(2\pi)$ jointly across all frequencies. $\square$

The results stated in Theorem 3.2 and its Corollary overlap to some extent with Theorem 2 of Li (2008). He uses the least absolute deviations estimator in the harmonic regression model

$$\hat{\beta}_n(\lambda) = \underset{(b_1,b_2)^\top \in \mathbb{R}^2}{\arg\min} \sum_{t=1}^n \rho_{0.5}\big(X_t - \cos(t\lambda)b_1 - \sin(t\lambda)b_2\big),$$

to define the *Laplace periodogram* $L_n(\lambda) = n|\hat{\beta}_n(\lambda)|^2/4$. In the special case that the time series of interest satisfies $X_t = \cos(t\lambda_0)\beta_1 + \sin(t\lambda_0)\beta_2 + \varepsilon_t$, where $\lambda_0$, $\beta_1$, and $\beta_2$ are unknown constants, this approach has the advantage that the maximizer of $L_n(\lambda)$ can be used as a robust estimator of $\lambda_0$, although Li provides only Monte Carlo evidence of this assertion. For general time series, he assumes that $X_t$ has median zero and a density $F'_X$ with $F'_X(0) > 0$, and that certain short-range dependence conditions are satisfied. The proofs of his Theorems 1 and 2 then yield an asymptotically linear representation for $\hat{\beta}_n(\lambda_n)$ that can be used to show

$$L_n(\lambda_n) = F'_X(0)^{-2} \Big|n^{-1/2} \sum_{t=1}^n \hat{V}_t(0.5) e^{-it\lambda_n}\Big|^2 + o_p(1).$$

The first term on the right-hand side is $2\pi/F'_X(0)^2$ times the quantile periodogram evaluated at the median. Hence, if the median of $X_t$ is indeed zero, the Laplace periodogram and the quantile periodogram at the median are asymptotically equivalent up to the unknown constant $2\pi/F'_X(0)^2$. Li (2012) extends his idea of harmonic median regression to quantile regression.

Using Li's (2008, 2012) periodograms instead of the quantile spectral methods introduced in my paper has the following disadvantages: (i) All of Li's asymptotic results depend on terms of the form $\tau(1-\tau)/F'_X(\xi(\tau))^2$ that in his case must be estimated to make inference about the dimensionless quantity $g_\tau(\lambda)$ even for simple tests such as in Example 3.4; my approach avoids this complication altogether. (ii) Li's methods require quantile regression at every frequency, whereas the quantile periodogram (2.3) can be computed easily with the Fast Fourier Transform. (iii) Li does not provide consistent estimators. For example, $L_n(\lambda)$ converges to a distribution with asymptotic mean $[2\pi/(4F'_X(0)^2)] \times g_{0.5}(\lambda)$, but Li does not establish that a smoothed version of $L_n(\lambda)$ converges in probability to this quantity. In contrast—as I will show now—the quantile periodogram can be smoothed by standard methods to provide uniformly consistent estimates of $g_\tau$.

Consistent estimation of the quantile spectrum requires weaker conditions than the con-struction of confidence intervals because much of the randomness introduced by replacing $r_\tau$ (as defined above (2.2)) with $\hat{r}_{n,\tau}$ is now controlled by the smoothing weight function $w$. Let $\varepsilon_0^*$ be an iid copy of $\varepsilon_0$ such that $X_t$ and $X_t^* := Y(\varepsilon_t, \ldots, \varepsilon_1, \varepsilon_0^*, \varepsilon_{-1}, \ldots)$ differ only through the input at time $t = 0$. I assume $X_t$ satisfies the following:

**Assumption C.** *For a given $\tau \in (0,1)$ and $\mathfrak{X}_\tau(\delta)$ as in Assumption A, there exists a $\delta > 0$ such that*
$$\sum_{t=0}^{\infty} \sup_{\xi \in \mathfrak{X}_\tau(\delta)} \|1\{X_t < \xi\} - 1\{X_t^* < \xi\}\| < \infty.$$

*Remarks.* 1. Assumption A implies Assumption C in view of the relation $\|1\{X_n' < \xi\} - 1\{X_n^* < \xi\}\| = \|1\{X_{n+1}' < \xi\} - 1\{X_{n+1} < \xi\}\|$; see the discussion below equation [13] of Wu (2005). For $\xi$ near $\xi_0(\tau)$, adding and subtracting $1\{X_n' < \xi\}$ and the triangle inequality then yield $\|1\{X_n < \xi\} - 1\{X_n^* < \xi\}\| = O(\sigma^n)$, which remains valid after taking suprema over $\mathfrak{X}_\tau(\delta)$.

2. A stationary stochastic process is usually called *short-range dependent* if its autocovariance function is summable. Since $X_t$ can have heavy tails, this definition no longer has the desired meaning. However, Remark 2.1 of Shao (2011a) can be used to show that $V_t(\tau)$ is short-range dependent because
$$\sum_{j \in \mathbb{Z}} |r_\tau(j)| \leq \left( \sum_{t=0}^{\infty} \|1\{X_t < \xi_0(\tau)\} - 1\{X_t^* < \xi_0(\tau)\}\| \right)^2 < \infty,$$

provided Assumptions A or C hold. This suggests that these assumptions should still be regarded as short-range dependence conditions on $X_t$. Heyde (2002) argues similarly to quantify the dependence of the increments of certain Gaussian processes.

Assumption C is easily verified in most cases via Proposition 3.1. However, more direct arguments can also be useful:

**Example 3.5** (Linear processes with Cauchy innovations)**.** Consider the linear process $X_t = \sum_{j=0}^{\infty} a_j \varepsilon_{t-j}$, where $(a_j)_{j \in \mathbb{N}}$ is a sequence of constants and $(\varepsilon_t)_{t \in \mathbb{Z}}$ is an sequence of iid copies of a standard Cauchy random variable. Without loss of generality, let $a_0 = 1$. Write $F_\varepsilon$ for the distribution function of $\varepsilon$; then $X_t$ has distribution function $F_X(x) = \mathbb{E} F_\varepsilon(x - \sum_{j=1}^{\infty} a_j \varepsilon_{t-j})$ and therefore also possesses a bounded density $F_X'$ by the Lebesgue Dominated Convergence Theorem. Furthermore, apply the point-wise inequality $|1\{X_n < \xi\} - 1\{X_n^* < \xi\}| \leq 1\{|X_n - \xi| < |X_n - X_n^*|\}$, then the Mean Value Theorem and $\mathbb{P}(|\varepsilon_0| + |\varepsilon_0^*| \geq x) \leq \mathbb{P}(|\varepsilon_0| \geq x/2) + \mathbb{P}(|\varepsilon_0^*| \geq x/2)$ for any fixed $x$ to see that

$$\begin{aligned}
\|1\{X_n < \xi\} - 1\{X_n^* < \xi\}\|^2 &\leq \mathbb{P}(|X_n - \xi| < |a_n||\varepsilon_0 - \varepsilon_0^*|) \\
&\leq \mathbb{P}(|X_0 - \xi| \leq |a_n|^{1/2}) + \mathbb{P}(|a_n||\varepsilon_0 - \varepsilon_0^*| \geq |a_n|^{1/2}) \\
&\leq 2|a_n|^{1/2} \sup_{x \in \mathbb{R}} F_X'(x) + 2\mathbb{P}(|\varepsilon_0| \geq |4a_n|^{-1/2}),
\end{aligned}$$

which is $O(|a_n|^{1/2})$ because the tail probability $\mathbb{P}(|\varepsilon_0| > x)$ of a Cauchy random variable is proportional to $x^{-1}$ as $x \to \infty$. Because these bounds hold uniformly in $\xi$, take square roots in the preceding display to conclude that Assumption C is satisfied if $\sum_{j=0}^{\infty} |a_j|^{1/4} < \infty$. The same type of reasoning can be used more generally when the innovations come from a smooth distribution whose tails behave algebraically. Proposition 3.1 does not apply because $(a_n)_{n \in \mathbb{N}}$ does not necessarily vanish at a geometric rate. □

Theorem 3.6 below establishes uniform consistency of the smoothed quantile periodogram under the condition that the bandwidth $B_n$ grows at a sufficiently slow rate. In particular, due to the uniformity, Theorem 3.6 allows for both fixed frequencies and sequences of frequencies such as the natural frequencies above.

**Theorem 3.6.** *If Assumptions B and C hold for some $\tau \in (0, 1)$, $w \in \mathcal{W}$, $B_n \to \infty$, and $B_n = o(\sqrt{n})$, then*

$$\hat{g}_{n,\tau}(\lambda) \xrightarrow{p} g_\tau(\lambda)$$

*uniformly in $\lambda \in (-\pi, \pi]$.*

*Remarks.* 1. The proof of Theorem 3.6 relies in part on recent results for classical spectral density estimates obtained by Liu and Wu (2010).

2. At fixed frequencies, kernel spectral density estimates of differentiable functions are often valid for bandwidths up to order $B_n = o(n)$; see, e.g., Andrews (1991) and Davidson and de Jong (2000). The stronger requirement $B_n = o(\sqrt{n})$ reflects that $\hat{V}_t(\tau)$ is not a smooth function of $\hat{\xi}_n(\tau)$. However, this requirement is not much of a restriction because, as Andrews notes, optimal bandwidths are typically of order less than $\sqrt{n}$.

3. If the quantile of interest $\xi_0(\tau)$ is assumed to be known, then Theorem 1 of Liu and Wu (2010) implies that Theorem 3.6 continues to hold when (i) $\xi_0(\tau)$ is used in the definition of $\hat{g}_{n,\tau}$ instead of $\hat{\xi}_n(\tau)$, (ii) Assumption B is replaced the condition that $F_X$ is continuous and increasing at $\xi_0(\tau)$, (iii) $\delta = 0$ in Assumption C, and (iv) $B_n = o(n)$.

4. The smoothed quantile periodogram at a known quantile $\xi_0(\tau)$ is just an ordinary smoothed periodogram of $V_t(\tau)$ and therefore optimality results from classical spectral analysis apply. In particular, the optimal lag window among the kernels in $\mathcal{W} \cap \{W \geq 0\}$ with respect to the relative mean-square error (MSE) criterion of Priestley (1962) is the quadratic-spectral (QS) window

$$w_{\text{QS}}(x) = \frac{25}{12\pi^2 x^2} \left( \frac{\sin(6\pi x/5)}{6\pi x/5} - \cos(6\pi x/5) \right).$$

The mean-square optimal bandwidth for the QS kernel is $B_n = O(n^{1/5})$, which can be established under additional dependence conditions; for example, Assumption A with $\delta = 0$ suffices. In the general case where $\xi_0(\tau)$ is estimated, a truncated MSE criterion as in Andrews (1991) could be used to limit the influence of $\hat{\xi}_n(\tau)$. However, his results rely crucially on second-order differentiability of the smoothed periodogram with respect to the estimated parameter. A fundamentally different approach is therefore likely to be needed, which I leave for future research.

A standard result in classical spectral analysis states that $\sqrt{n/B_n}$ times the centered smoothed periodogram converges to a normal distribution. Establishing a similar result for the smoothed quantile periodogram is complicated by the dependence of the spectral density estimate on the sample quantile $\hat{\xi}_n(\tau)$. One way of simplifying the dependence structure is to estimate the quantile and the associated smoothed quantile periodogram separately using approximately one half of the data each. Let $l_n = (\log n)^p$ for some $p > 1$ and $m_n = \lfloor n/2 \rfloor$; the constant $p$ can be chosen arbitrarily close to 1. Denote the new quantile estimate based on the first $m_n - l_n$ observations by $\tilde{\xi}_n(\tau) := \hat{\xi}_{m_n - l_n}(\tau)$. The observations $m_n + 1 \leq t \leq n$ are then used to compute the smoothed quantile periodogram

$$\tilde{g}_{n,\tau}(\lambda) = \frac{1}{2\pi m_n} \sum_{|j| < m_n} w(j/B_n) \cos(j\lambda) \sum_{t=|j|+m_n+1}^{n} V_t\bigl(\tau, \tilde{\xi}_n(\tau)\bigr) V_{t-|j|}\bigl(\tau, \tilde{\xi}_n(\tau)\bigr).$$

The remaining $l_n$ observations are discarded. They facilitate the approximation of the data by $l_n$-dependent variables $\tilde{X}_t = \mathbb{E}(X_t \mid \varepsilon_t, \varepsilon_{t-1}, \ldots, \varepsilon_{t-l_n+1})$. Assume the following:

**Assumption D.** *There is some $n^*$ such that for all $n > n^*$, $F_{\tilde{X}}(x) := \mathbb{P}(\tilde{X}_0 \leq x)$ is Lipschitz continuous in a neighborhood of $\xi_0(\tau)$ and $\mathbb{E}|X_n - X'_n| = O(\varrho^n)$ for some $\varrho \in (0, 1)$.*

This assumption is similar to the conditions of Proposition 3.1; however, instead of the Lipschitz continuity of the limit $F_X$ of $F_{\tilde{X}}$, it requires $F_{\tilde{X}}$ to be eventually Lipschitz. Assumption D and some additional conditions on the lag window imply the asymptotic normality of the modified smoothed quantile periodogram:

**Theorem 3.7.** *Suppose Assumptions B and D hold for some $\tau \in (0, 1)$, $w$ is even and Lipschitz continuous with support $[-1, 1]$, $w(0) = 1$, $\lim_{x \to 0}(1 - w(x))/|x|^3 < \infty$, $B_n \to \infty$, $B_n = o(n^{1/4})$, $n = o(B_n^7)$. Then*

$$\sqrt{m_n/B_n}\bigl(\tilde{g}_{n,\tau}(\lambda) - g_\tau(\lambda)\bigr) \rightsquigarrow \mathrm{N}\bigl(0, \sigma^2(\lambda)\bigr),$$

*where $\sigma^2(\lambda) = (1 + h(2\lambda))g_\tau(\lambda)^2 \int_{-1}^{1} w(x)^2\, dx$, and $h(\lambda) = 1$ if $\lambda = 2\pi k$ for some $k \in \mathbb{Z}$ and $0$ otherwise.*

I investigate the finite sample properties of the smoothed quantile periodogram and confidence intervals based on the quantile periodogram in a small simulation study in section 5. The next section discusses the use of integrated quantile periodograms to test for uninformative quantile spectra.

## 4. Testing for Flatness of a Quantile Spectrum

In this section I provide two Cramér-von Mises tests (Procedures 4.4 and 4.6 below) for the null hypothesis that the $\tau$-th quantile spectrum is flat, i.e., $g_\tau(\lambda) \equiv \tau(1-\tau)/(2\pi)$, against the alternative that the $\tau$-th quantile spectrum is informative.

If the distribution function of $X_t$ is continuous and increasing at $\xi_0(\tau)$, then $r_\tau(0) =$

$\tau(1-\tau)$ and the null and alternative hypotheses can be stated more precisely as

$$\text{H}_0\colon r_\tau(j) = 0 \text{ for all } j > 0 \quad \text{and} \quad \text{H}_1\colon r_\tau(j) \neq 0 \text{ for some } j > 0.$$

Provided that $\sum_{j \in \mathbb{Z}} r_\tau(j)$ converges absolutely, the $\tau$-th quantile spectrum is symmetric about zero. One way to test for the null hypothesis is therefore to check if the sample equivalent of

$$\int_0^\lambda g_\tau(u)\,du - \int_0^\lambda \frac{r_\tau(0)}{2\pi}\,du = \sum_{j>0} r_\tau(j)\psi_j(\lambda), \tag{4.1}$$

where $\psi_j(\lambda) := \sin(j\lambda)/(\pi\lambda)$, is near zero for all $\lambda \in \Pi := [0, \pi]$.

The quantity in the preceding display is best understood as an function in $L_2(\Pi)$, the set of Lebesgue-measurable functions $f\colon \Pi \to \mathbb{R}$ with $\int_\Pi f(\lambda)^2\,d\lambda < \infty$. Under the equivalence relation "$f \equiv g$ if and only if $f = g$ Lebesgue-almost everywhere," $L_2(\Pi)$ is a proper Hilbert space with inner product $\langle f, g\rangle_\Pi := \int_\Pi f(\lambda)g(\lambda)\,d\lambda$ for $f, g \in L_2(\Pi)$ and norm $\|f\|_\Pi := \sqrt{\langle f, f\rangle_\Pi}$. Since $\|\psi_j\|_\Pi^2 = 1/(2\pi j^2)$ for all $j \in \mathbb{Z} \setminus \{0\}$, (4.1) indeed lies in $L_2(\Pi)$ and satisfies

$$\left\|\sum_{j>0} r_\tau(j)\psi_j\right\|_\Pi^2 = \sum_{j>0} r_\tau(j)^2 \|\psi_j\|_\Pi^2.$$

Here we need the fact that $\langle \psi_j, \psi_k\rangle_\Pi = 0$ for all $j \neq k$. Now replace $r_\tau(j)$ by $\hat{r}_{n,\tau}(j)$ and rescale to obtain the Cramér-von Mises statistic

$$CM_{n,\tau} := \left\|\sqrt{n}\sum_{j=1}^{n-1} \hat{r}_{n,\tau}(j)\psi_j\right\|_\Pi^2 = \frac{n}{2\pi}\sum_{j=1}^{n-1}\left(\frac{\hat{r}_{n,\tau}(j)}{j}\right)^2$$

based on the random process $S_{n,\tau}(\lambda) := \sqrt{n}\sum_{j=1}^{n-1} \hat{r}_{n,\tau}(j)\psi_j(\lambda)$ in $L_2(\Pi)$. No smoothing weight function and bandwidth is needed because the integral in (4.1) already acts as a smoothing operator. The scaling factor $\sqrt{n}$ in $S_{n,\tau}$ is included because $\sqrt{n}\hat{r}_{n,\tau}(j)$ can be expected to have an asymptotically normal distribution for each $j > 0$ under the null hypothesis. When viewed as a random function on $L_2(\Pi)$, the process $S_{n,\tau}(\lambda)$ should then converge in distribution to a mean-zero Gaussian process $S_\tau(\lambda)$ with covariances

$$\mathbb{E} S_\tau(\lambda) S_\tau(\lambda') = \sum_{j>0}\sum_{k>0}\sum_{l\in\mathbb{Z}} \text{Cov}\big(V_0(\tau)V_j(\tau), V_{j-l}(\tau)V_{j-l-k}(\tau)\big)\psi_j(\lambda)\psi_k(\lambda'), \tag{4.2}$$

$\lambda, \lambda' \in \Pi$, so that $CM_{n,\tau} \rightsquigarrow \|S_\tau\|_\Pi^2$ by the Continuous Mapping Theorem (see, e.g., Theorem 18.11 of van der Vaart, 1998, p. 259).

As the following theorem shows, this convergence indeed occurs if the conditions of the null hypothesis are strengthened slightly: Suppose that under $\text{H}_0$, for a given $\tau \in (0, 1)$ there is a $\delta > 0$ such that

$$\mathbb{P}(X_0 < \xi, X_j < \xi') = \mathbb{P}(X_0 < \xi)\mathbb{P}(X_0 < \xi') \text{ for all } j > 0 \text{ and all } \xi, \xi' \in \mathfrak{X}_\tau(\delta), \quad (4.3)$$

where $\mathfrak{X}_\tau(\delta) = \{\xi \in \mathbb{R} : |\xi_0(\tau) - \xi| \leq \delta\}$ as before. The role of this condition is discussed in detail in the remarks and examples below.

**Theorem 4.1.** *Suppose Assumptions A and B hold for some $\tau \in (0, 1)$.*
  (i) *If $H_0$ is satisfied in the sense of (4.3), then $CM_{n,\tau} \rightsquigarrow \|S_\tau\|_\Pi^2$, and*
  (ii) *if $H_1$ is satisfied, then $\mathbb{P}(CM_{n,\tau} > B) \to 1$ for every $B \in \mathbb{R}$.*

*Remarks.* 1. For the proof of the theorem, I show stochastic equicontinuity of the empirical process $(n-j)^{-1/2} \sum_{j=1}^{n-j} [V_t(\tau, \xi)V_{t+j}(\tau, \xi) - \mathbb{E}V_0(\tau, \xi)V_j(\tau, \xi)]$ indexed by $\xi$ under Assumptions A and B for each fixed $j$. Condition (4.3) is used to control the behavior of $\hat{r}_{n,\tau}(j)$ for large $j$ and $n$. These two results then allow me to apply a general result on Cramér-von Mises tests for spectral densities given in Shao (2011a).

2. Condition (4.3) imposes slightly more on the dependence structure of $V_t(\tau)$ than the white noise assumption $H_0$ (i.e., $\delta = 0$). However, since $\delta$ can be chosen to be as small as desired, it is much less restrictive than requiring that $(X_t)_{t\in\mathbb{Z}}$ be pairwise independent ($\delta = \infty$) or even iid, which is frequently imposed when testing for white noise; see, e.g., Milhøj (1981) and Hong (1996).

**Example 4.2** (Stochastic volatility, continued)**.** Recall that $F_\varepsilon$ is the distribution function of $\varepsilon$. The stochastic volatility process in Example 2.1 has a flat $\tau_0$-th quantile spectrum but fails to satisfy (4.3) because

$$\mathbb{P}(X_0 < \xi, X_j < \xi') = \mathbb{E}1\left\{\varepsilon_0 < \frac{\xi - \xi_0(\tau_0)}{v(\varepsilon_{-1}, \ldots)}\right\} F_\varepsilon\left(\frac{\xi' - \xi_0(\tau_0)}{v(\varepsilon_{j-1}, \ldots)}\right)$$

can generally not be simplified further due to the lagged innovations contained in the volatility process $v$. Thus, Theorem 4.1 does not apply. However, the test procedure from Example 3.4 can still be used in this case to test for flatness of the $\tau_0$-th quantile spectrum, for if $g_{\tau_0}(\lambda_0) = \tau_0(1-\tau_0)/(2\pi)$ is rejected at some frequency $\lambda_0$, then $H_1$ must be true. Linton and Whang (2007) investigate the stochastic volatility model of Example 2.1 with the sample *quantilogram*, defined as $\hat{r}_{n,\tau}(j)/\hat{r}_{n,\tau}(0)$, for a fixed, finite number of lags $j = 1, 2, \ldots$. From their results it can be seen that the failure of (4.3) for the stochastic volatility model manifests itself in terms of a non-vanishing drift term in $\sqrt{n}\hat{r}_{n,\tau}(j)$ due to the estimation of $\xi_0(\tau)$. A Cramér-von Mises test requires control of these drifts for large $j$ and $n$; this is nontrivial and left for future work. □

**Example 4.3** (QAR, continued)**.** The QAR process in Example 2.2 possesses a flat $\tau_0$-th quantile spectrum and has the property (4.3) if there exists a neighborhood $\mathcal{T}$ of $\tau_0$ such that $\beta_1(\tau) = \beta_2(\tau) = 0$ for all $\tau \in \mathcal{T}$: In this case, the conditional quantile function, defined as the solution $\xi(\tau \mid \mathcal{F}_{t-1})$ of $\mathbb{P}(X_t \leq \xi \mid \mathcal{F}_{t-1}) = \tau$, is given by $\xi(\tau \mid \mathcal{F}_{t-1}) = \beta_0(\tau) + \beta_1(\tau)X_{t-1} + \beta_2(\tau)X_{t-2} = \beta_0(\tau)$ almost surely for all $\tau \in \mathcal{T}$ by monotonicity. Take expectations to deduce that

$$\tau = \mathbb{P}\big(X_t \leq \xi(\tau \mid \mathcal{F}_{t-1}) \mid \mathcal{F}_{t-1}\big) = \mathbb{P}\big(X_t \leq \beta_0(\tau)\big) = \mathbb{P}\big(X_0 \leq \xi_0(\tau)\big),$$

almost surely for all $\tau \in \mathcal{T}$ and therefore $\xi(\tau \mid \mathcal{F}_t) = \xi_0(\tau)$ almost surely on $\tau \in \mathcal{T}$. Conclude that for any $\tau, \tau' \in \mathcal{T}$,

$$\begin{aligned}\mathbb{P}\big(X_0 < \xi_0(\tau), X_j < \xi_0(\tau')\big) &= \mathbb{E}1\{X_0 < \xi_0(\tau)\}\mathbb{P}\big(X_j < \xi_0(\tau') \mid \mathcal{F}_{j-1}\big) \\ &= \mathbb{P}\big(X_0 < \xi_0(\tau)\big)\mathbb{P}\big(X_0 < \xi_0(\tau')\big).\end{aligned}$$

Now (4.3) follows because as long as $F_X$ is continuous and increasing in a neighborhood of $\xi_0(\tau_0)$, there is a $\delta > 0$ such that for every $\xi, \xi' \in \mathcal{X}_{\tau_0}(\delta)$, there are $\tau, \tau' \in \mathcal{T}$ such that $\xi = \xi_0(\tau)$ and $\xi' = \xi_0(\tau')$. The assertion in Example 2.2 about the flatness of the $\tau_0$-th quantile spectrum is obtained by letting $\mathcal{T} = \{\tau_0\}$. □

The main difficulty with applying Theorem 4.1 in practice is the unknown covariance function (4.2) of the limiting process $S_\tau$. In standard spectral analysis, this has led researchers to assume that $X_t$ is iid normal under the null hypotheses of white noise (Durbin, 1967, is an important early reference) in order to avoid having to estimate the covariance function of a Gaussian process. In sharp contrast, in quantile spectral analysis the assumption that $X_t$ is iid is already enough to construct a test for flatness without imposing a distributional assumption: In large samples $\hat{V}_t(\tau)$ is close to $V_t(\tau) = \tau - 1\{X_t < \xi_0(\tau)\}$ in probability, but $1\{X_t < \xi_0(\tau)\}$ is a Bernoulli random variable with success probability $\tau$ as long as $F_X$ is continuous and increasing at $\xi_0(\tau)$. Hence, if $X_t$ is indeed iid and $J_1, J_2, \ldots, J_n$ are independent Bernoulli$(\tau)$ variables, then

$$\tilde{CM}_{n,\tau} := \frac{1}{2\pi n}\sum_{j=1}^{n-1} j^{-2}\left(\sum_{t=1+j}^{n} V_t(\tau)V_{t-j}(\tau)\right)^2 \quad \text{and}$$

$$CM'_{n,\tau} := \frac{1}{2\pi n}\sum_{j=1}^{n-1} j^{-2}\left(\sum_{t=1+j}^{n} (\tau - J_t)(\tau - J_{t-j})\right)^2$$

have the same distribution. Because $CM_{n,\tau} = \tilde{CM}_{n,\tau} + o_p(1)$ under the conditions of Theorem 4.1(i), this distributional equivalence leads to a simple, distribution-free Monte Carlo test. I prove its consistency in Corollary 4.5 below.

**Procedure 4.4** (Monte Carlo test for flatness). 1. Draw $n$ iid copies $J_1, J_2, \ldots, J_n$ of a Bernoulli$(\tau)$ random variable.
2. Compute $CM'_{n,\tau}$ with the variables from step 1.
3. Repeat steps 1 and 2 $R$ times. Reject $H_0$ in favor of $H_1$ if $CM_{n,\tau}$ is larger than $c_{n,\tau}(1-\alpha)$, the $1-\alpha$ empirical quantile of the $R$ realizations of $CM'_{n,\tau}$.

*Remark.* Exploiting the distribution-free character of sign or quantile crossing indicators has a long history in statistics and econometrics; see, e.g., Walsh (1960). More recently,

Chernozhukov, Hansen, and Jansson (2009) use it to construct finite sample confidence intervals for quantile regression estimators.

By choosing the number of Monte Carlo repetitions $R$ large enough, the quantiles of the null distribution of $\tilde{CM}_{n,\tau}$ can be approximated with arbitrary precision. I therefore let $R \to \infty$ and define the quantiles of the simulated distribution directly as $c_{n,\tau}(1 - \alpha) := \inf\{x \in \mathbb{R} : \mathbb{P}(\tilde{CM}_{n,\tau} > x) \leq \alpha\}$. The large sample properties of Procedure 4.4 can now be stated as follows:

**Corollary 4.5.** *Suppose Assumption B holds for some $\tau \in (0,1)$ and let $\alpha \in (0,1)$.*
  (i) *If $(X_t)_{t \in \mathbb{Z}}$ is an iid sequence, then $\mathbb{P}(CM_{n,\tau} > c_{n,\tau}(1 - \alpha)) \to \alpha$, and*
  (ii) *if instead Assumption A and $H_1$ are satisfied, then $\mathbb{P}(CM_{n,\tau} > c_{n,\tau}(1 - \alpha)) \to 1$.*

*Remark.* If $\xi_0(\tau)$ is known, then the test in Procedure 4.4 has level $\alpha$ even in finite samples provided that $\tilde{CM}_{n,\tau}$ is used in step 3 instead of $CM_{n,\tau}$.

In cases where it does not seem reasonable to assume that $X_t$ is iid under the null hypothesis, the block-wise wild bootstrap of Shao (2011a) should be used instead. This bootstrap is a modification of the standard wild bootstrap (Liu, 1988; Mammen, 1992). It perturbs whole blocks of observations with iid copies of a random variable $\eta$ that is independent of the data and satisfies $\mathbb{E}\eta = 0$, $\mathbb{E}\eta^2 = 1$, and $\mathbb{E}\eta^4 < \infty$. Since the blocks grow with the sample size, this eventually captures enough of the dependence structure to provide critical values for the null distribution under the more general condition (4.3).

**Procedure 4.6** (Shao's block-wise wild bootstrap). 1. Choose a block length $b_n \leq n$ and the corresponding number of blocks $L_n = n/b_n$, taken to be an integer for convenience. For each $s = 1, \ldots L_n$ define a block $\mathcal{B}_s = \{(s-1)b_n + 1, \ldots, sb_n\}$.
  2. Draw $L_n$ iid copies $\eta_1, \eta_2, \ldots, \eta_{L_n}$ of $\eta$. For each $t = 1, \ldots, n$, define $\omega_t = \sum_{s=1}^{L_n} \eta_s 1\{t \in \mathcal{B}_s\}$ so that $\omega_t$ takes on the value $\eta_s$ if $t$ lies in the $s$-th block.
  3. Compute $\hat{r}^*_{n,\tau}(j) := n^{-1} \sum_{t=j+1}^n [\hat{V}_t(\tau)\hat{V}_{t-j}(\tau) - \hat{r}_{n,\tau}(j)]\omega_t$ and calculate the bootstrap statistic
$$CM^*_{n,\tau} := \frac{n}{2\pi} \sum_{j=1}^{n-1} \left(\frac{\hat{r}^*_{n,\tau}(j)}{j}\right)^2.$$
  4. Repeat steps 2 and 3 $R$ times. Reject $H_0$ in favor of $H_1$ if $CM_{n,\tau}$ is larger than $c^*_{n,\tau}(1 - \alpha)$, the $1 - \alpha$ empirical quantile of the $R$ realizations of $CM^*_{n,\tau}$.

*Remark.* The recommended choice for $\eta$ in practice is a Rademacher variable that takes on the value 1 with probability 1/2 and the value $-1$ with probability 1/2. Distributions other than the Rademacher distribution can be used for $\eta$, in particular if $\hat{V}_t(\tau)\hat{V}_{t-j}(\tau)$ has a skewed distribution, but there is no evidence that they would lead to better inference; see Davidson, Monticini, and Peel (2007) for a discussion of this point for the standard wild bootstrap.

As before, I take $R$ to be large and define the quantiles of the bootstrap distribution conditional on the sample $\mathcal{S}_n$ as $c^*_{n,\tau}(1 - \alpha) = \inf\{x \in \mathbb{R} : \mathbb{P}(CM^*_{n,\tau} \leq x \mid \mathcal{S}_n) \geq 1 - \alpha\}$. Procedure 4.6 then has the following asymptotic properties:

**Theorem 4.7.** *Suppose Assumptions A and B hold for some $\tau \in (0,1)$. Let $\alpha \in (0,1)$,*

$b_n \to \infty$ and $b_n/n \to 0$.
  (i) If $H_0$ is satisfied in the sense of (4.3), then $\mathbb{P}(CM_{n,\tau} > c^*_{n,\tau}(1-\alpha)) \to \alpha$, and
  (ii) if $H_1$ is satisfied, then $\mathbb{P}(CM_{n,\tau} > c^*_{n,\tau}(1-\alpha)) \to 1$.

*Remark.* If $\xi_0(\tau)$ is known, then Theorems 4.1 and 4.7 remain valid without condition (4.3) as long as $\tilde{CM}_{n,\tau}$ is used in place of $CM_{n,\tau}$.

The next section investigates the finite sample properties of the two Cramér-von Mises tests, the quantile periodogram, and the smoothed quantile periodogram in a Monte Carlo study and provides an empirical application.

## 5. Numerical Results

In this section I present a sequence of examples to illustrate quantile spectral methods in the context of some familiar time series models and macroeconomic data, and compare the results to those obtained from traditional spectral analysis.

**Example 5.1** (AR(2) with spectral peak). Let $(\varepsilon_t)_{t \in \mathbb{Z}}$ be iid copies of an $N(0,1)$ variable with distribution function $\Phi$. Li (2008) investigates the frequency domain properties of a stationary AR(2) process of the form

$$X_t = \beta_1 X_{t-1} + \beta_2 X_{t-2} + \varepsilon_t, \qquad \beta_1 = 2 \times 0.95 \cos(2\pi \times 0.22), \beta_2 = -0.95^2. \qquad (5.1)$$

Shao and Wu's (2007) Theorem 5.2 implies that $X_t$ is GMC for all $\alpha > 0$. Since $X_t$ is also normally distributed, Proposition 3.1 applies and consequently Assumptions A and B hold. To study the finite sample properties of classical and quantile spectral estimates for sample sizes $n \in \{300, 600, 900\}$ in this model, I generated 10,000 realizations of the process of size $400 + n$ for each $n$ and then discarded the first 400 observations. Each realization was initialized by independent standard normal random variables. The solid black line in the left panel of Figure 1 plots a QS-smoothed periodogram of $X_t$, i.e.,

$$\hat{f}_{n,X}(\lambda) = \frac{1}{2\pi} \sum_{|j|<n} w_{\text{QS}}(j/B_n) \hat{\gamma}_{n,X}(j) \cos(j\lambda),$$

where $\hat{\gamma}_{n,X}(j) := n^{-1} \sum_{t=|j|+1}^{n} (X_t - \bar{X}_n)(X_{t-|j|} - \bar{X}_n)$ and $\bar{X}_n := n^{-1} \sum_{t=1}^{n} X_t$, of one such realization with $n = 300$ and $B_n = 13n^{1/5} \approx 40.68$. The process (5.1) has little noise and a single pronounced peak at $2\pi \times 0.22$ in its spectral density, shown as the dotted line in the left panel of Figure 1. The smoothed periodogram therefore does not have much difficulty identifying the peak, although its size is underestimated slightly due to the smoothing. The shaded area in the left panel shows 95% asymptotic point-wise confidence bands based on the *periodogram* of $X_t$, defined as

$$I_{n,X}(\lambda) = \frac{1}{2\pi} \sum_{|j|<n} \hat{\gamma}_{n,X}(j) \cos(j\lambda).$$

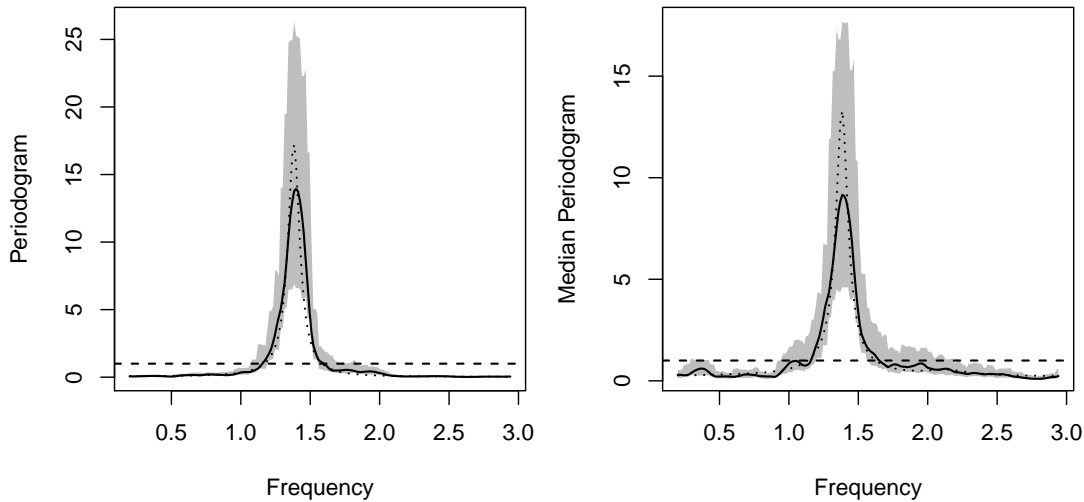

Figure 1: Left panel: spectral density (dotted line) of $X_t$ as in (5.1), QS-smoothed periodogram (solid) of a realization with $n = 300$ and $B_n = 13n^{1/5} \approx 40.68$, chi-squared point-wise 95% confidence bands (shaded grey) with $k = 4$, and $\hat{\gamma}_{n,X}(0)/(2\pi)$ (dashed). Right panel: median spectrum (dotted) of $X_t$, QS-smoothed median periodogram (solid), chi-squared point-wise 95% confidence bands (shaded grey), and $0.5(1 - 0.5)/(2\pi)$ (dashed). Both panels use the same data, $B_n$, and $k$, and are normalized by $\hat{\gamma}_{n,X}(0)/(2\pi)$ (left) and $0.5(1 - 0.5)/(2\pi)$ (right).

The point-wise confidence bands were computed by averaging over $2k + 1$ periodogram coordinates at natural frequencies in the same way as in Corollary 3.3, but with $Q_{n,\tau}$ replaced by $I_{n,X}$. Here and in all plots below, I used $k = 4$. The dashed line in the left panel plots $\hat{\gamma}_{n,X}(0)/(2\pi)$, i.e., the usual estimate of $f_X$ if the spectrum were known to be flat. It provides a natural point of comparison for the other quantities; in particular, it can be seen from the left panel that the peak at $2\pi \times 0.22$ is significantly different from a flat spectrum at the 5% level.

The right panel of Figure 1 analyzes the same data with quantile spectral methods. The black line is the QS-smoothed median (i.e., 0.5-th quantile) periodogram and the shaded area graphs 95% point-wise confidence bands computed as described in Corollary 3.3. Here I used the same values for $B_n$ and $k$ as in the left panel. The dashed line is $0.5(1-0.5)/(2\pi)$, i.e., the median spectrum under the hypothesis that it is flat. The dotted line shows the median spectrum $g_{0.5}$, which can be calculated exactly from equation (6) in Li (2008). The smoothed median periodogram clearly identifies the peak, although the estimate of the actual size of the peak is slightly worse than the one obtained in the left panel. However, the median spectrum is completely contained inside the confidence bands and the peak at $2\pi \times 0.22$ differs significantly from a flat median spectrum at the 5% level.

For both panels the choice of $B_n$ and $k$ matters, with lower values of $B_n$ and higher values of $k$ leading to smoother—but not necessarily better—estimates: Figure 2 shows the mean integrated square error (MISE) of the QS-smoothed periodogram (left panel) and the QS-smoothed median periodogram (right) estimated from the 10,000 realizations as a function of $B_n/n^{1/5}$. Here the behavior of both methods is quite similar and the MISEs attain their minimum at $B_n/n^{1/5} \approx 13$ for each $n \in \{300, 600, 900\}$, which provides evidence that the optimal growth rate $B_n = O(n^{1/5})$ for the QS-smoothed periodogram is also a good choice for QS-smoothed quantile periodograms. Further, Table 1 shows the empirical frequency of the event that the 95% confidence interval at $\lambda \in \{\pi \times 0.22, 2\pi \times 0.22, 3\pi \times 0.22\}$ covered the spectrum and median spectrum, respectively, in the experiments for $k \in \{2, 4, 6\}$ and $n$ as before. The confidence intervals constructed from the periodogram and the median periodogram behaved very similar at the three frequencies and covered the population value in nearly 95% of all cases unless $n$ was small and $k$ was large. For these values both methods had low coverage frequencies. □

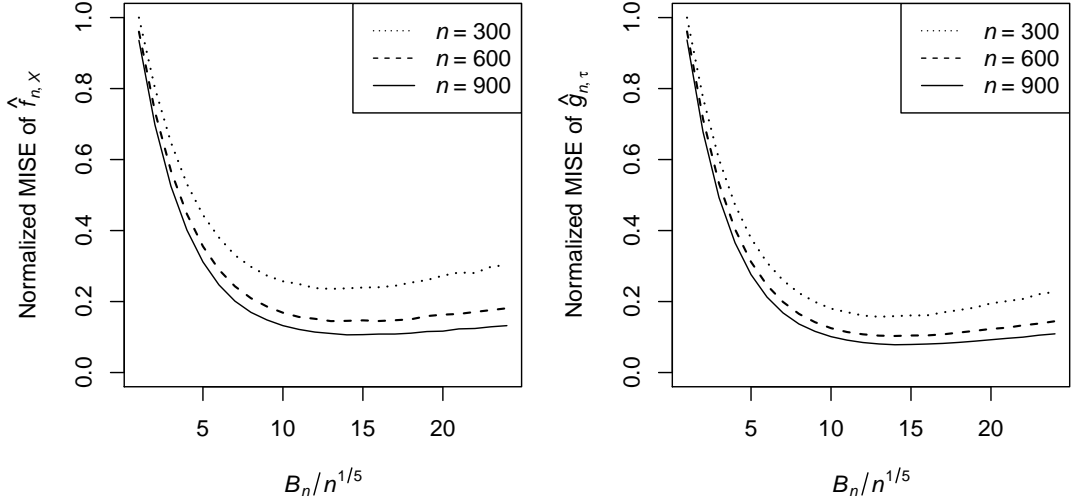

Figure 2: Empirical MISE of the QS-smoothed periodogram (left panel) and the QS-smoothed median periodogram (right) for three different sample sizes as a function of $B_n/n^{1/5}$. Both panels were separately normalized by the respective joint maximum of the three curves.

Robust estimators (in the sense of Huber and Ronchetti, 2009, p. 5) exhibit *stability*, i.e., small deviations from the model assumptions should have small effects on the performance of the estimator, and high *breakdown resistance*, i.e., larger deviations should not cause catastrophic results. The following two examples illustrate that classical spectral estimates are not robust to outliers in the data, whereas quantile spectral estimators provide reliable results in such situations.

**Example 5.2** (Stability of quantile spectral estimators). Suppose that each observation in a realization of the AR(2) process from Example 5.1 has a probability $p$ of being contaminated by an additional additive error component. For this I drew iid Bernoulli($p$) variables $J_1, \ldots, J_n$ and iid central Student t($\nu$) variables $T_1, \ldots, T_n$ to generate the observed samples as $\mathcal{S}_n = \{X_t + J_t T_t : t = 1, \ldots, n\}$, where the $X_1, \ldots, X_n$ were taken from Example

Table 1: Finite-sample coverage frequencies of an asymptotic 95% confidence interval (CI) for the spectrum and median spectrum of the process in Example 5.1 at $\lambda \in \{\pi \times 0.22, 2\pi \times 0.22, 3\pi \times 0.22\}$ as a function of $n$ and $k$.

|  |  | Periodogram CI | | | Median Periodogram CI | | |
| --- | --- | --- | --- | --- | --- | --- | --- |
| $n$ | $k$ | $\pi \times .22$ | $2\pi \times .22$ | $3\pi \times .22$ | $\pi \times .22$ | $2\pi \times .22$ | $3\pi \times .22$ |
| 300 | 2 | 0.940 | 0.931 | 0.937 | 0.937 | 0.982 | 0.961 |
|  | 4 | 0.936 | 0.676 | 0.931 | 0.924 | 0.907 | 0.974 |
|  | 6 | 0.921 | 0.249 | 0.909 | 0.913 | 0.178 | 0.979 |
| 600 | 2 | 0.943 | 0.951 | 0.948 | 0.942 | 0.982 | 0.956 |
|  | 4 | 0.946 | 0.915 | 0.947 | 0.938 | 0.980 | 0.962 |
|  | 6 | 0.944 | 0.774 | 0.941 | 0.926 | 0.930 | 0.965 |
| 900 | 2 | 0.950 | 0.951 | 0.948 | 0.948 | 0.974 | 0.956 |
|  | 4 | 0.949 | 0.941 | 0.946 | 0.940 | 0.980 | 0.959 |
|  | 6 | 0.948 | 0.904 | 0.947 | 0.934 | 0.971 | 0.964 |

5.1. The spectral density of the corresponding process $(X_t + J_t T_t)_{t \in \mathbb{Z}}$ is

$$f_{X+JT}(\lambda) = f_X(\lambda) + \frac{p}{2\pi} \frac{\nu}{\nu - 2},$$

which, for any given $p$, can be made as large as desired by choosing $\nu > 2$ sufficiently close to 2 without violating the assumptions of classical spectral theory. Figure 3 plots $f_{X+JT}(\lambda)$ for $p = 0.15$ and $\nu = 2.001$ as a dotted line in the left panel; the median spectrum (dotted, right) needed no adjustment because it is invariant under such contamination. The other quantities are the same as in Figure 1 and the same 300 observations were used, but 46 of these were contaminated. The smoothed periodogram retains the spectral shape and has a significant spike at $2\pi \times 0.22$, but grossly underestimates the location of the spectrum. Moreover, the confidence bands no longer contain the spectrum at any frequency. In sharp contrast, the smoothed median periodogram is barely affected by the contamination and the confidence bands cover the median spectrum at almost all frequencies. The hypothesis that $g_{0.5}(2\pi \times 0.22) = 0.5(1-0.5)/(2\pi)$ can also be clearly rejected.

The odd behavior of the classical spectral density estimates in Figure 3 is likely due to the imprecisely estimated auto-covariances of the contaminated process. As Basraka, Davis, and Mikosch (2002) point out, for near-infinite variance time series the convergence rate of sample auto-covariances to their population equivalent is much slower than $n^{-1/2}$. Since periodograms are weighted sums of sample auto-covariances, they can be expected to inherit this lack of precision. In contrast, the sample auto-covariances of $\hat{V}_t(\tau)$ can be shown to converge at rate $n^{-1/2}$ as long as Assumption A and a slightly strengthened version of Assumption B hold.

I repeated the experiment from Table 1 with the contaminated data. The estimated coverage probabilities for the confidence intervals constructed from the periodogram and the median periodogram are shown in Table 2. As can be seen, the presence of outliers had little effect on the performance of the quantile spectral estimates. In sharp contrast, the coverage probability for the classical spectrum was almost zero in most cases and 0.245 in the best scenario ($k = 2$, $n = 900$). □

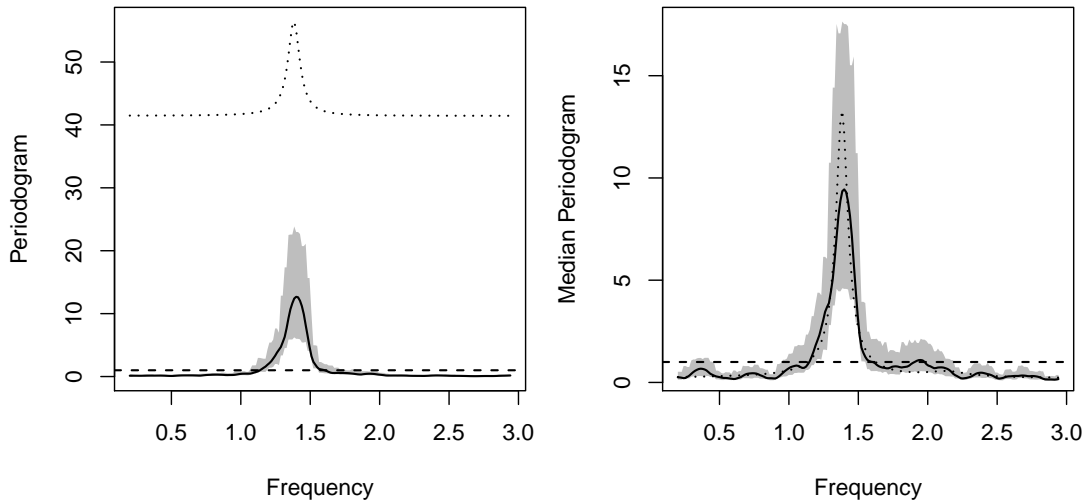

Figure 3: Left panel: spectral density (dotted line) of the process in Example 5.2, QS-smoothed periodogram (solid) of a realization with $n = 300$ and $B_n = 13n^{1/5} \approx 40.68$, chi-squared point-wise 95% confidence bands (shaded grey) with $k = 4$, and $\hat{\gamma}_{n,X}(0)/(2\pi)$ (dashed). Right panel: median spectrum (dotted), QS-smoothed median periodogram (solid), chi-squared point-wise 95% confidence bands (shaded grey), and $0.5(1 - 0.5)/(2\pi)$ (dashed). Both panels use the same data, $B_n$, and $k$, and are normalized by $\hat{\gamma}_{n,X}(0)/(2\pi)$ (left) and $0.5(1 - 0.5)/(2\pi)$ (right).

**Example 5.3** (Breakdown resistance of quantile spectral estimators)**.** Now suppose instead that each observation from Example 5.1 has a 15 percent chance of being contaminated by one of the iid Cauchy$(0, 1)$ variables $C_1, \ldots, C_n$. The observed samples then were $\mathcal{S}_n = \{X_t + J_t C_t : t = 1, \ldots, n\}$ with the $X_1, \ldots, X_n$ as before. Since these outliers do not have a well defined mean, the spectral density of the corresponding contaminated process no longer exists. Spectral analysis by ordinary methods broke down completely when 46 of the 300 observations used for Figure 1 were contaminated: The smoothed periodogram in Figure 4 no longer has the expected spectral shape and fails to give any indication of a periodicity present in the data. A comparison of the confidence bands to the estimate of $\gamma_X(0)/(2\pi)$ now provides overwhelming evidence for the false hypothesis that the process is white noise. In sharp contrast, the median spectrum is unaffected by the contamination and the smoothed median periodogram significantly identifies the periodicity. In addition, the confidence bands remain essentially unchanged from Example 5.2, which is also confirmed by the coverage probability estimates of the confidence intervals constructed from median periodograms provided in Table 3. Here the estimates were nearly identical to the ones presented in Table 2 for the median spectrum. Corresponding estimates for the classical spectrum cannot be computed because it is unbounded at all frequencies. □

Table 2: Finite-sample coverage frequencies of an asymptotic 95% confidence interval (CI) for the spectrum and the median spectrum of the process in Examples 5.2 as a function of $n$ and $k$.

| $n$ | $k$ | Periodogram CI | | | Median Periodogram CI | | |
|---|---|---|---|---|---|---|---|
| | | $\pi \times .22$ | $2\pi \times .22$ | $3\pi \times .22$ | $\pi \times .22$ | $2\pi \times .22$ | $3\pi \times .22$ |
| 300 | 2 | 0.001 | 0.109 | 0.001 | 0.918 | 0.976 | 0.958 |
| | 4 | 0.001 | 0.001 | 0.001 | 0.900 | 0.858 | 0.965 |
| | 6 | 0.001 | 0.001 | 0.001 | 0.886 | 0.089 | 0.972 |
| 600 | 2 | 0.001 | 0.208 | 0.001 | 0.928 | 0.976 | 0.952 |
| | 4 | 0.001 | 0.006 | 0.001 | 0.914 | 0.966 | 0.958 |
| | 6 | 0.001 | 0.001 | 0.001 | 0.901 | 0.890 | 0.963 |
| 900 | 2 | 0.001 | 0.245 | 0.001 | 0.932 | 0.974 | 0.948 |
| | 4 | 0.001 | 0.013 | 0.001 | 0.924 | 0.973 | 0.957 |
| | 6 | 0.001 | 0.002 | 0.001 | 0.912 | 0.953 | 0.957 |

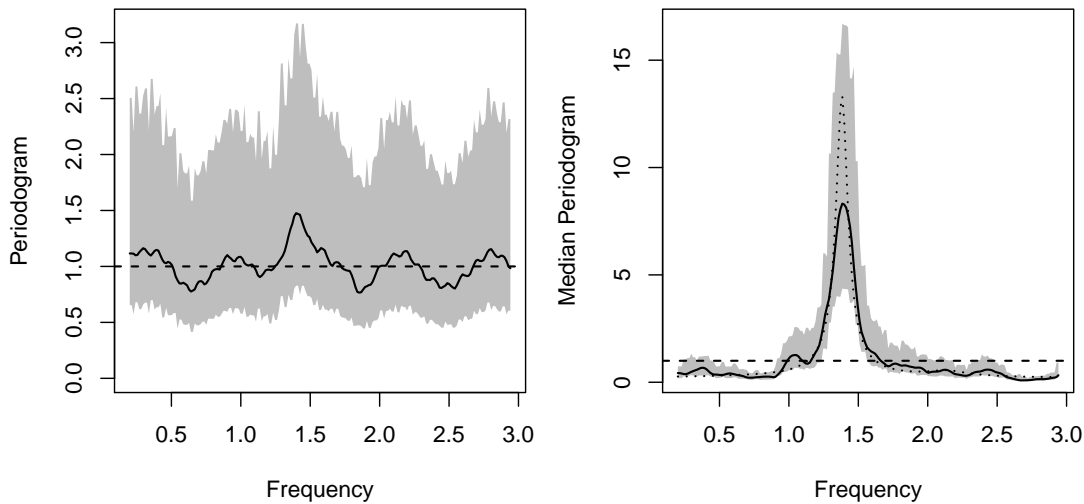

Figure 4: Left panel: QS-smoothed periodogram (solid black) of a realization of the process in Example 5.3 with $n = 300$ and $B_n = 13n^{1/5} \approx 40.68$, chi-squared point-wise 95% confidence bands (shaded grey) with $k = 4$, and $\hat{\gamma}_{n,X}(0)/(2\pi)$ (dashed). The spectral density does not exist. Right panel: median spectrum (dotted), QS-smoothed median periodogram (solid), chi-squared point-wise 95% confidence bands (shaded grey), and $0.5(1 - 0.5)/(2\pi)$ (dashed). Both panels use the same data, $B_n$, and $k$, and are normalized by $\hat{\gamma}_{n,X}(0)/(2\pi)$ (left) and $0.5(1 - 0.5)/(2\pi)$ (right).

Table 3: Finite-sample coverage frequencies of an asymptotic 95% confidence interval for the median spectrum of the process in Example 5.3 as a function of $n$ and $k$.

|  |  | Median Periodogram | | |
| --- | --- | --- | --- | --- |
| $n$ | $k$ | $\pi \times .22$ | $2\pi \times .22$ | $3\pi \times .22$ |
| 300 | 2 | 0.904 | 0.971 | 0.957 |
|  | 4 | 0.889 | 0.827 | 0.960 |
|  | 6 | 0.861 | 0.059 | 0.966 |
| 600 | 2 | 0.915 | 0.973 | 0.948 |
|  | 4 | 0.901 | 0.961 | 0.952 |
|  | 6 | 0.884 | 0.872 | 0.955 |
| 900 | 2 | 0.918 | 0.970 | 0.951 |
|  | 4 | 0.905 | 0.966 | 0.952 |
|  | 6 | 0.883 | 0.940 | 0.951 |

For the next Monte Carlo exercise, I return to the stochastic volatility model from Example 2.1 to illustrate that even if the classical spectrum shows no sign of periodicity, almost all quantiles of the distribution can be crossed in a periodic manner.

**Example 5.4** (Stochastic volatility, continued). Take $(\varepsilon_t)_{t\in\mathbb{Z}}$ to be iid copies of an $N(0, \theta^2)$ variable and let $u_t = \log v(\varepsilon_{t-1}, \varepsilon_{t-2}, \ldots)$ be the stationary solution of the process $u_t = \beta_1 u_{t-1} + \beta_2 u_{t-2} + \varepsilon_{t-1}$ with $\beta_1, \beta_2$ as in (5.1). Then $e^{u_t}$ is log-normally distributed and $X_t = \varepsilon_t v(\varepsilon_{t-1}, \varepsilon_{t-2}, \ldots) = \varepsilon_t e^{u_t}$ has median zero. To show that $X_t$ is GMC, apply the Mean Value Theorem and the Cauchy-Schwarz inequality to obtain the bound $\|X_n - X'_n\|_\alpha \leq \|\varepsilon_n\|_\alpha \|e^{\bar{u}_n}\|_{2\alpha} \|u_n - u'_n\|_{2\alpha}$, where $u'_n$ is $u_n$ with $(\varepsilon_0, \varepsilon_{-1}, \ldots)$ replaced by $(\varepsilon^*_0, \varepsilon^*_{-1}, \ldots)$ and $\bar{u}_n$ lies on the line segment joining $u_n$ and $u'_n$. By monotonicity of the exponential function and the Minkowski inequality, we have

$$\|e^{\bar{u}_n}\|_{\max\{1,2\alpha\}} \leq \|\max\{e^{u_n}, e^{-u_n}, e^{u'_n}, e^{-u'_n}\}\|_{\max\{1,2\alpha\}} \leq 4\|e^{u_n}\|_{\max\{1,2\alpha\}} < \infty$$

because the four terms inside the maximum have the same log-normal distribution. If needed, the Loève $c_r$ inequality provides a similar bound for the case $0 < 2\alpha < 1$. The GMC property then follows since $u_t$ is GMC by Theorem 5.2 of Shao and Wu (2007). The distribution function of $X_t$ is given by $F_X(x) = \mathbb{E}\Phi(x/(e^{u_t}\theta))$, which can be seen to have a bounded density with the help of the Lebesgue Dominated Convergence Theorem. Therefore, Assumptions A and B again hold.

The top two panels of Figure 5 graph the same spectral estimates as in Figure 1 for $n = 600$ observations of the stochastic volatility model with $\theta = 1$. The spectrum (not shown to prevent clutter) and the median spectrum (identical to the dashed line in the top right panel) of the model are flat, which is also correctly identified at almost all frequencies by both point-wise confidence bands. The bottom two panels show the smoothed quantile periodograms (black lines) and point-wise confidence bands (shaded grey) at $\tau = 0.25$ (left) and $\tau = 0.75$ (right) computed from the same data. In both panels, the estimated quantile spectra show a considerable spike that is significantly different from a flat $\tau$-th quantile spectrum at frequency $2\pi \times 0.22$, thereby providing evidence of a dependence structure that is not present in the mean and auto-covariance of the process. In addition, the estimated quantile spectra at $\tau = 0.25$ and $\tau = 0.75$ have a large, but less informative spike at frequency zero (not shown to enhance readability) that corresponds to the long-run variation discovered at quantiles other than the median. Since the quantile spectra of the process do not possess a closed-form expression for $\tau \neq 0.5$, I instead also plot smoothed quantile periodograms of $n = 10^6$ observations at $\tau = 0.25$ (left) and $\tau = 0.75$ (right) as dotted lines in the bottom panels to illustrate how much of the spectral shape is already recovered in a sample with 600 observations. Indeed, although the estimates from the smaller sample are more volatile, the size and shape of the peaks at $2\pi \times 0.22$ are nearly identical for the two sample sizes.

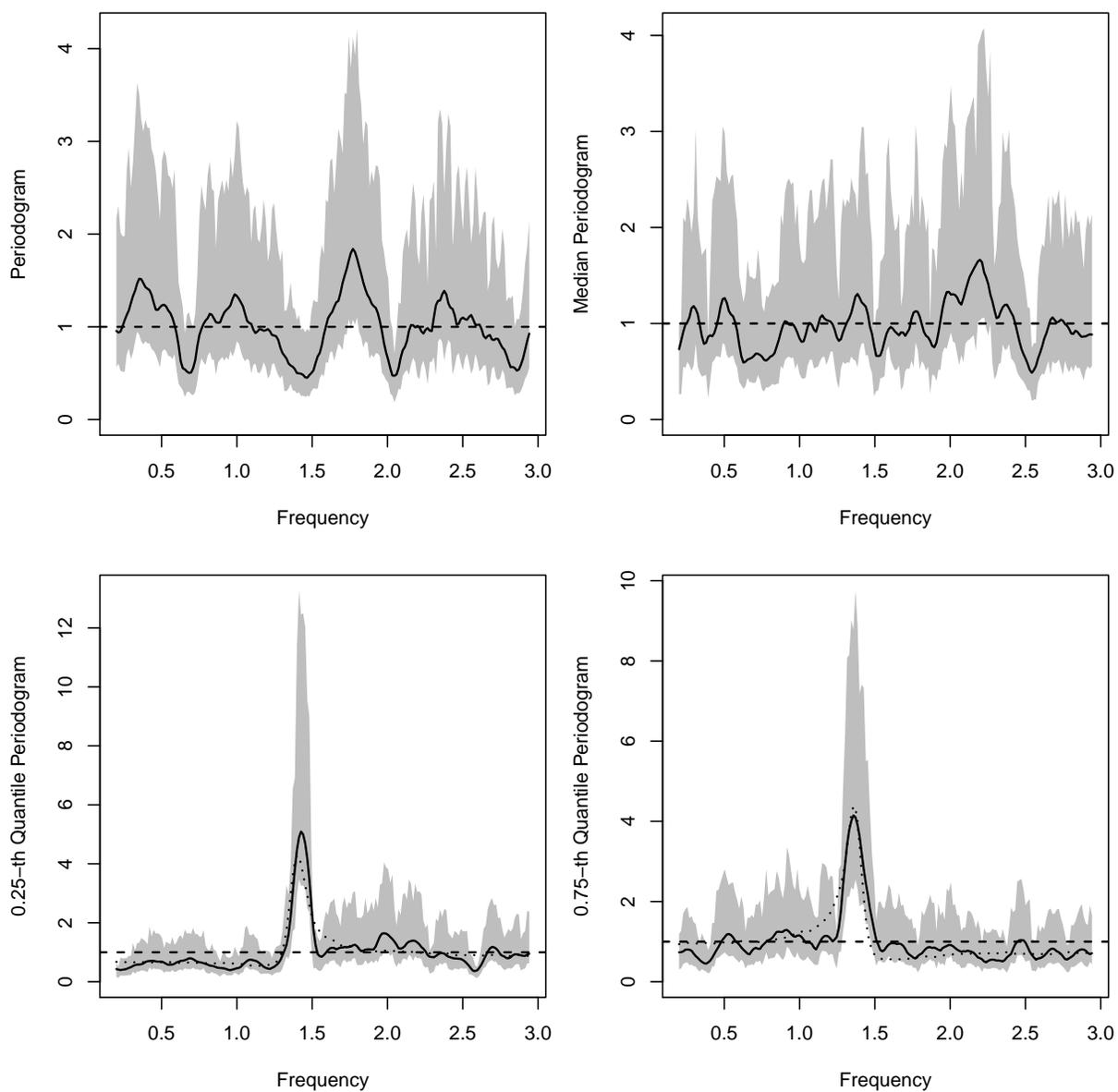

Figure 5: Top left panel: QS-smoothed periodogram (solid) of a realization of the process in Example 5.4 with $n = 600$ and $B_n = 13n^{1/5} \approx 46.73$, chi-squared point-wise 95% confidence bands (shaded grey) with $k = 4$, and $\hat{\gamma}_{n,X}(0)/(2\pi)$ (dashed). Other panels: QS-smoothed $\tau$-th quantile periodogram (solid), chi-squared point-wise 95% confidence bands (shaded grey), and $\tau(1-\tau)/(2\pi)$ (dashed) for $\tau = 0.5$ (top right), 0.25 (bottom left), and 0.75 (bottom right). All panels use the same data, $B_n$, and $k$. The top left panel is normalized by $\hat{\gamma}_{n,X}(0)/(2\pi)$. The other panels are normalized by $\tau(1-\tau)/(2\pi)$. The bottom two panels also show QS-smoothed $\tau$-th quantile periodograms (dotted) with $n = 10^6$ for $\tau = 0.25$ (left) and 0.75 (right). Frequencies near zero are not shown to enhance readability.

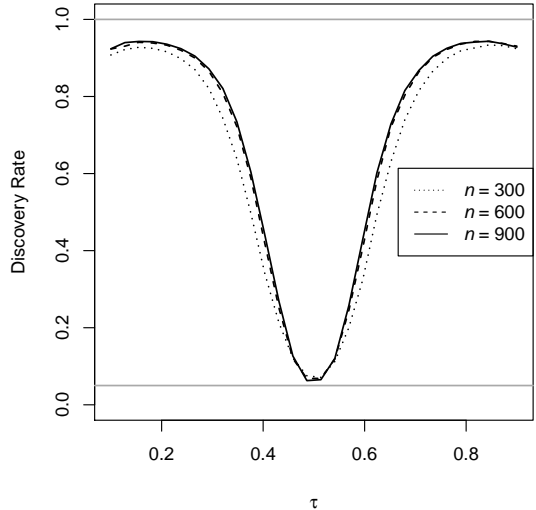

Figure 6: Empirical size and power of a test for a cycle with frequency $2\pi \times 0.22$ in the stochastic volatility model of Example 5.4 as a function of $\tau$. Nominal size at $\tau = 0.5$ is 0.05 (lower grey line).

To evaluate how reliably the quantile spectral estimates discover the cycle at frequency $2\pi \times 0.22$, I recorded the relative number of the test decisions in favor of the hypothesis $H_0\colon g_\tau(2\pi \times 0.22) = \tau(1-\tau)/(2\pi)$ in 10,000 realizations of the stochastic volatility model using a 95% confidence interval with $k = 4$. The results are shown in Figure 6 for different sample sizes as a function of $\tau \in (0,1)$. At $\tau = 0.5$, the null hypothesis is true and the tests almost attained the 5% level (lower grey line) for the three sample sizes. At the other quantiles, the null hypothesis is false, which was also correctly recognized at all sample sizes as long as a quantile not too close to $\tau = 0.5$ was chosen. In particular, near the quartiles the power of the tests was about 0.9.

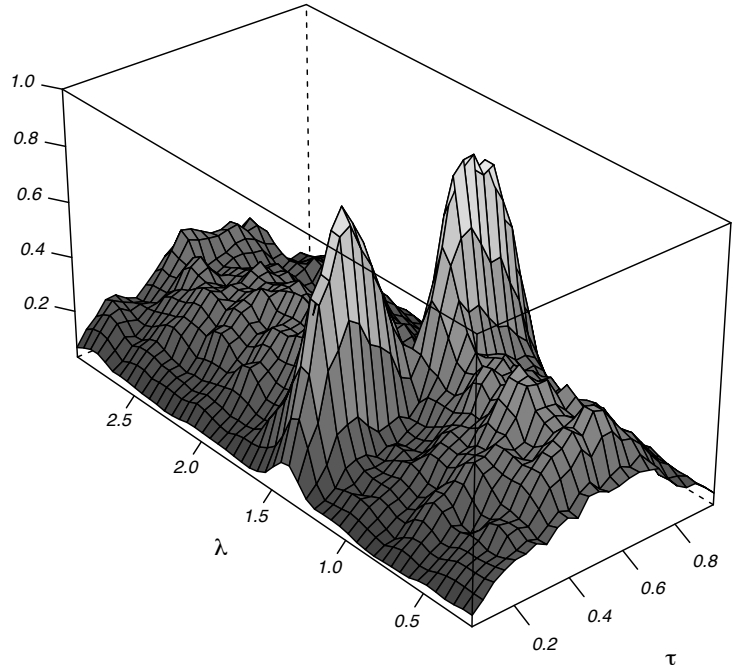

Figure 7: QS-smoothed quantile periodogram across all quantiles of a realization of the process in Example 5.4 with $n = 900$ and $B_n = 8n^{1/5} \approx 31.18$, normalized by the joint maximum of all quantile periodograms. Frequencies near zero are not shown for better readability.

The additional information obtained from quantile spectral analysis can also be seen in Figure 7, where I graph the QS-smoothed quantile periodogram as a function of both $\lambda$ and $\tau$. Here I chose $n = 900$ and $B_n = 8n^{1/5} \approx 31.18$ for a smoother appearance of the plot. The two humps in the figure make it clear that most of the dependence structure is in fact present near the lower and upper quartiles of the process, whereas working with the mean or median provides no insight in this case. □

The following examples illustrate the size and power of the two Cramér-von Mises tests introduced in section 4.

**Example 5.5** (QAR(2) and Procedure 4.4)**.** Table 4 shows the empirical rejection frequency of the null hypothesis of a flat $\tau$-th quantile spectrum as a function of $n \in \{100, 200, 300\}$ and $\tau \in \{0.1, 0.5, 0.9\}$ in a variety of settings. For each entry, I recorded the test decision of Procedure 4.4 in 10,000 realizations by comparing the test statistics to 5% critical values obtained from $10^6$ simulations each. The first column of the "Size" portion provides the rejection frequencies when the data were iid $\chi_3^2$ variables. In this case, the null hypothesis is true at all quantiles. The test behaved mildly conservatively for $\tau = 0.1$ in smaller samples, but was close to the level of the test at other quantiles and samples sizes. In samples larger than 300 (not reported), the test was essentially exact at all quantiles. I also experimented with other distributions, including normal, Student t(2), and standard Cauchy variables, but found that they had little impact on the results.

Table 4: Rejection frequencies of the null hypothesis for the Monte Carlo Cramér-von Mises test (Procedure 4.4) at the 5% level.

|  |  | Size | | | Power | | |
| --- | --- | --- | --- | --- | --- | --- | --- |
| $n$ | $\tau$ | $\chi_3^2$ | Ex. 5.4 | QAR | Ex. 5.1 | Ex. 5.4 | QAR |
| 100 | 0.1 | 0.022 | – | 0.024 | 0.093 | 0.007 | – |
|  | 0.5 | 0.053 | 0.068 | – | 0.999 | – | 0.999 |
|  | 0.9 | 0.037 | – | – | 0.169 | 0.332 | 0.993 |
| 200 | 0.1 | 0.019 | – | 0.021 | 0.405 | 0.043 | – |
|  | 0.5 | 0.052 | 0.076 | – | 1.000 | – | 1.000 |
|  | 0.9 | 0.046 | – | – | 0.504 | 0.468 | 1.000 |
| 300 | 0.1 | 0.048 | – | 0.029 | 0.795 | 0.188 | – |
|  | 0.5 | 0.052 | 0.080 | – | 1.000 | – | 1.000 |
|  | 0.9 | 0.050 | – | – | 0.875 | 0.724 | 1.000 |

The first column of the "Power" portion shows the relative number of rejections when the data-generating process was the AR(2) from Example 5.1. Here the null hypothesis is false at all quantiles, which was also reliably identified at the median at all samples. However, at the outer quantiles the spectral peak is smaller and therefore larger samples were needed to detect its presence. The results for the contaminated processes from Examples 5.2 and 5.3 are not shown because they were virtually identical.

The second "Size" and "Power" columns give the rejection frequencies for the stochastic volatility model from Example 5.4. The null hypothesis is true at $\tau = 0.5$, but the process is not covered by the assumptions underlying the Monte Carlo test because the stochastic volatility model is not iid, which resulted in a mild over-rejection at all sample sizes. At the other quantiles, the process satisfies $H_1$ and the test has power against this alternative by Corollary 4.5(ii). The power of the test increased sharply with the sample size for $\tau = 0.9$, whereas for $\tau = 0.1$ the increase was considerably slower. Some intuition for this result can be gathered from Figure 5, where the estimated quantile spectrum in the lower quantiles can be seen to have a long stretch on which it is close to the hypothetical quantile spectrum implied by the null hypothesis. In contrast, this stretch is somewhat shorter in the upper quantiles. Moreover, as shown in Figure 7, the setup for the test is quite demanding because the spectral peak near the extremes of the distribution is small. Larger samples (not reported) yielded better results, with the power being nearly one at all quantiles for $n = 600$.

The third columns of the "Size" and "Power" portions show the relative number of rejections of the hypothesis of a flat $\tau$-th quantile spectrum for realizations of the QAR(2) process (see Example 2.2)

$$X_t = \underbrace{4 + \Phi^{-1}(\varepsilon_t)}_{\beta_0(\varepsilon_t)} + \underbrace{0.8 \times 1\{\varepsilon_t > 0.2\}}_{\beta_1(\varepsilon_t)} X_{t-1} + \underbrace{0.6 \times 1\{\varepsilon_t > 0.6\}}_{\beta_2(\varepsilon_t)} X_{t-2} \qquad (5.2)$$

where, as before, $(\varepsilon)_{t\in\mathbb{Z}}$ is a sequence of iid copies of a Uniform$(0,1)$ variable. By Theorem 5.1 of Shao and Wu (2007), this recursion admits a stationary solution of the form (3.1) and satisfies the GMC property. Further, the marginal distribution function of $X_t$ can be seen to possess a bounded Lebesgue density from the properties of truncated normal variables and dominated convergence. If $(X_t)_{t\in\mathbb{Z}}$ is positive, the right-hand side of (5.2) is guaranteed to be increasing in $\varepsilon_t$ conditional on $X_{t-1}, X_{t-2}$ and the model in the preceding display is indeed a proper QAR model. Since the process has a very small probability of generating a negative observation, I therefore considered only positive realizations of (5.2) in order to enforce well-behaved sample paths.

The QAR process satisfies the null hypothesis of a flat quantile spectrum for $\tau \in (0, 0.2]$ and the alternative at the other quantiles. In particular, it behaves like a stationary QAR(1) on $\tau \in (0.2, 0.6]$ that exhibits enough mean reversion to regulate the explosive behavior of the process on $\tau \in (0.6, 1)$. This dependence structure induces an asymmetric spectral shape across quantiles, with spectral peaks of different sizes at frequency zero in the middle to upper quantiles. The QS-smoothed quantile periodogram of a realization with $n = 900$ plotted in Figure 8 illustrates this shape. As can be seen from Table 4, the Monte Carlo Cramér-von Mises test very reliably detected the presence of the alternative hypothesis at $\tau = 0.5$ and $0.9$ even for $n = 100$. At $\tau = 0.1$ the null hypothesis is true and, although Procedure 4.4 does not apply because the observations are not iid, the test was only mildly conservative. □

**Example 5.6** (QAR(2) and Procedure 4.6)**.** I repeated the experiments outlined in the previous example with the wild bootstrap test described in Procedure 4.6. I experimented with the block size $b_n$, but found that the results were not overly sensitive to this choice as long as the blocks were not too large. I therefore settled for block sizes near $\sqrt{n}/2$ and used $b_n = 5$, $8$, and $10$ for $n = 100$, $200$, and $300$, respectively, although other choices are clearly possible; see Shao (2011a) for a thorough discussion.

The results are shown in Table 5. The important difference to the preceding example is that the QAR(2) model (5.2) is now fully covered by the assumptions of the test; see Theorem 4.7. This is also reflected in the test for a flat quantile spectrum of the QAR process at $\tau = 0.1$, which was nearly exact for $n = 300$. The other results in the "Size" portion of the table were similar to the ones given in Table 4 for the Monte Carlo test. The power of the bootstrap test was also comparable to the other test, but neither of the tests dominated the other: For the AR model both test behaved similarly, for the stochastic volatility model the bootstrap test showed a more balanced performance, and for the QAR model the Monte Carlo test was more powerful at the outer quantiles. □

**Example 5.7** (Building permits data)**.** Finally, to illustrate what kind of insights quantile spectral analysis of actual economic data can provide, I consider the series "New Privately Owned Housing Units Authorized by Building Permits in Permit-Issuing Places" from the US Census.[1] The data consist of 634 (seasonally unadjusted) monthly observations from January 1959 to October 2011 of the total number of permits from permit-issuing places in the United States that report to the Census. Such a permit is typically issued by a town or a county and enables an individual to begin construction on a new housing unit.

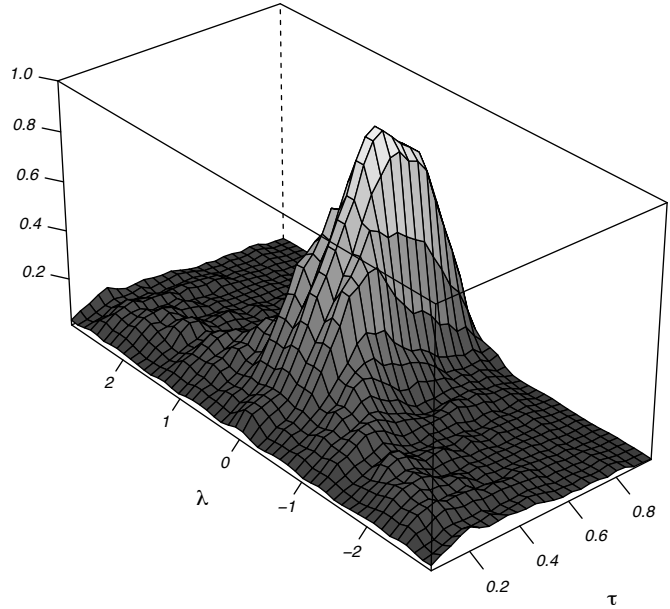

Figure 8: QS-smoothed quantile periodogram across all quantiles of a realization of the QAR(2) process in Example 5.5 with $n = 900$ and $B_n = 8n^{1/5} \approx 31.18$, normalized by the joint maximum of all quantile periodograms.

Table 5: Rejection frequencies of the null hypothesis for the bootstrap Cramér-von Mises test (Procedure 4.6) at the 5% level. The block sizes are $b_n = 5, 8$, and $10$ for $n = 100, 200$, and $300$, respectively. I used the warp-speed method of Giacomini et al. (2007) to estimate size and power of the bootstrap test; this method considerably sped up the simulations because only one bootstrap replication per Monte Carlo replication was needed.

|     |      | Size |         |       | Power   |         |       |
| --- | ---  | ---  | ---     | ---   | ---     | ---     | ---   |
| $n$ | $\tau$ | $\chi_3^2$ | Ex. 5.4 | QAR | Ex. 5.1 | Ex. 5.4 | QAR |
| 100 | 0.1  | 0.027 | –      | 0.026 | 0.113   | 0.056   | –     |
|     | 0.5  | 0.055 | 0.094  | –     | 1.000   | –       | 0.999 |
|     | 0.9  | 0.029 | –      | –     | 0.170   | 0.110   | 0.374 |
| 200 | 0.1  | 0.031 | –      | 0.030 | 0.430   | 0.339   | –     |
|     | 0.5  | 0.058 | 0.083  | –     | 1.000   | –       | 1.000 |
|     | 0.9  | 0.049 | –      | –     | 0.486   | 0.422   | 0.635 |
| 300 | 0.1  | 0.050 | –      | 0.051 | 0.754   | 0.550   | –     |
|     | 0.5  | 0.056 | 0.090  | –     | 1.000   | –       | 1.000 |
|     | 0.9  | 0.052 | –      | –     | 0.780   | 0.567   | 0.820 |

Figure 9 graphs this time series in the frequency domain: The smoothed periodogram (solid line, top left) and smoothed median periodogram (solid, top right) behave similarly and have their largest peaks at frequencies 0.045 and 0.039, respectively, which translates into an estimated business cycle length of 11.58 years when measured by the smoothed periodogram and 13.39 years when measured by the smoothed median periodogram. Both lines also have peaks of similar size at the yearly ($2\pi/12 \approx 0.52$) and half-yearly ($2\pi/6 \approx 1.05$) frequencies, which provides evidence of considerable seasonality in the data. However, as illustrated by the smoothed 0.10-th quantile periodogram (solid, bottom left) and 0.90-th quantile periodogram (solid, bottom right), these seasonal cycles do not appear uniformly across the distribution of the data. At the 0.90-th quantile the yearly and—to some extent—the half-yearly cycles are still present, but at the 0.10-th quantile this seasonality vanishes completely. The smoothed 0.10-th quantile periodogram also has some smaller peaks between 0.2 and 0.5, but comparison to the confidence intervals (shaded grey) shows that these peaks are not significantly different from a straight line. All graphs have in common, however, that the business cycle explains most of the cyclical variation, which indicates that the influence of seasonal patterns disappears during economic troughs. □

## 6. Conclusion

In this paper I introduced quantile spectral densities that summarize the cyclical behavior of time series across their whole distribution by analyzing periodicities in quantile crossings. I discussed robust spectral estimation and inference in situations where the dependence structure of a time series is not accurately captured by the auto-covariance function, in particular when the time series under consideration is uncorrelated or heavy-tailed. I established the statistical properties of quantile spectral estimators in a large class of nonlinear time series models and discussed inference both at fixed and across all frequencies. Monte Carlo experiments and an empirical example showed that quantile spectral estimates are similar to regular spectral density estimates in both shape and interpretation when standard conditions are satisfied, but can still reliably identify dependence structures when these conditions fail to hold.

---

[1] I downloaded the data from `http://www.census.gov/const/permits_cust.xls` on November 30, 2011.

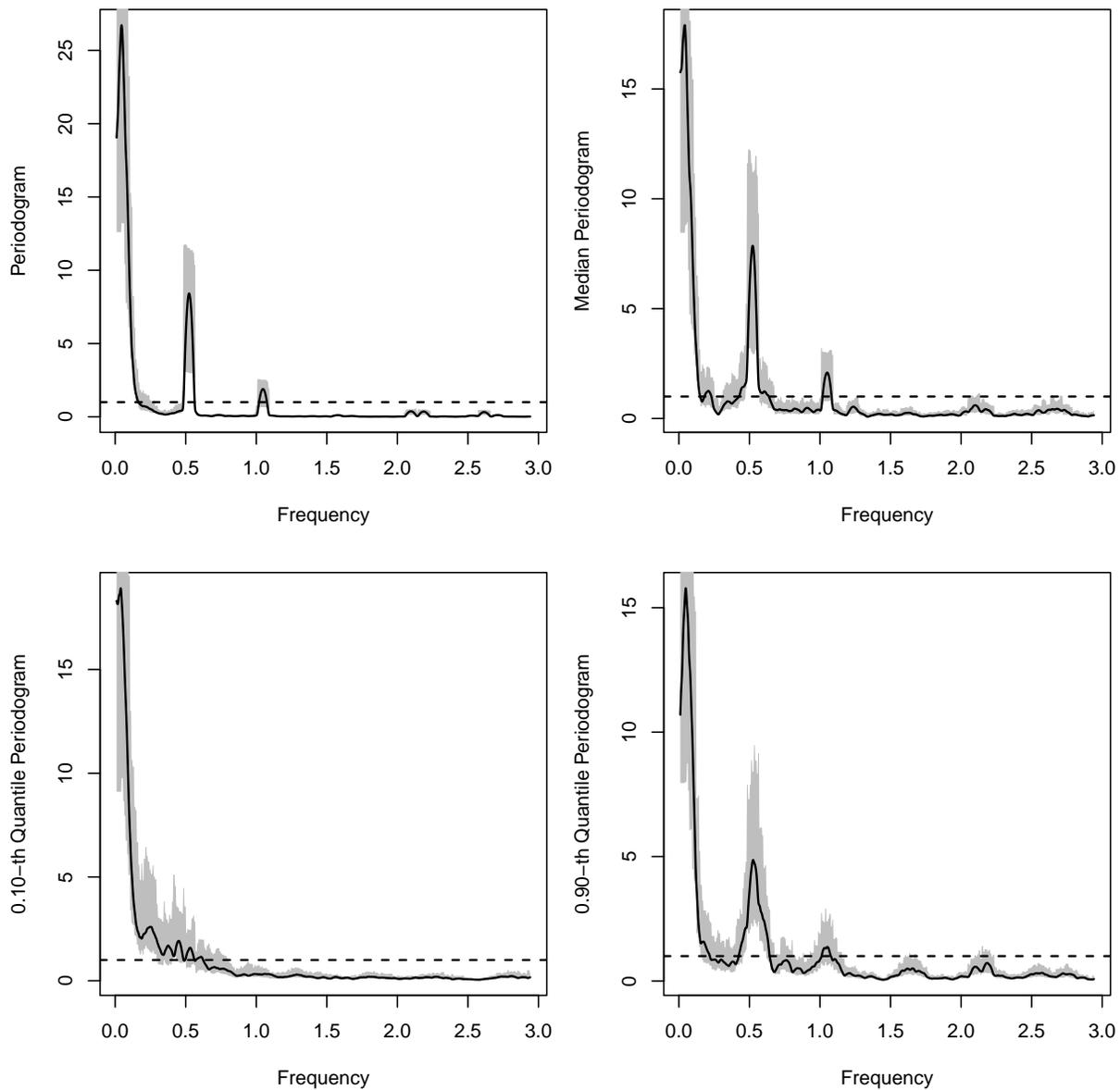

Figure 9: Top left panel: QS-smoothed periodogram (solid) of building permit data with $B_n = 28n^{1/5} \approx 101.76$ to prevent smaller peaks from being smoothed out, chi-squared point-wise 95% confidence bands (shaded grey, not shown completely to prevent clutter) with $k = 4$, and $\hat{\gamma}_{n,X}(0)/(2\pi)$ (dashed). Other panels: QS-smoothed $\tau$-th quantile periodogram (solid), chi-squared point-wise 95% confidence bands (shaded grey), and $\tau(1-\tau)/(2\pi)$ (dashed) for $\tau = 0.5$ (top right), 0.10 (bottom left), and 0.90 (bottom right). All panels use the same data, $B_n$, and $k$. The top left panel is normalized by $\hat{\gamma}_{n,X}(0)/(2\pi)$. The other panels are normalized by $\tau(1-\tau)/(2\pi)$.

# Appendix

Throughout the Appendix, $\mathbb{P}^*$ and $\mathbb{E}^*$ respectively denote outer probability and outer expectation (see, e.g., van der Vaart, 1998, p. 258). Probability and expectation conditional on the observed sample $\mathcal{S}_n$ is abbreviated by $\hat{\mathbb{P}}(\cdot) := \mathbb{P}(\cdot \mid \mathcal{S}_n)$ and $\hat{\mathbb{E}}(\cdot) := \mathbb{E}(\cdot \mid \mathcal{S}_n)$.

## A. Proofs

*Proof of Proposition 3.1.* By assumption, we can find a $\delta' > 0$ such that $F_X$ is Lipschitz on $(\xi_0(\tau) - \delta', \xi_0(\tau) + \delta')$. Choose a large enough $N \in \mathbb{N}$ such that $\delta := \delta' - \varrho^{N/(1+\alpha)} > 0$ and pick any $\xi \in \mathcal{X}_\tau(\delta)$; then, for all $n \geq N$, apply the pointwise bound $|1\{X_n < \xi\} - 1\{X_n' < \xi\}| \leq 1\{|X_n - \xi| < |X_n - X_n'|\}$, the Markov inequality, and the GMC property to see that

$$\|1\{X_n < \xi\} - 1\{X_n' < \xi\}\|^2 \leq \mathbb{P}(|X_n - \xi| < |X_n - X_n'|)$$
$$\leq \mathbb{P}(|X_n - \xi| < \varrho^{n\alpha/(1+\alpha)}) + \mathbb{E}|X_n - X_n'|^\alpha \varrho^{-n\alpha^2/(1+\alpha)}$$
$$\leq M\varrho^{n\alpha/(1+\alpha)}$$

for a large enough absolute constant $M$. This constant can be enlarged slightly to ensure that the inequality also holds for the remaining $n < N$. With $\sigma := \varrho^{\alpha/(2+2\alpha)}$, take square-roots on both sides and suprema over $\mathcal{X}_\tau(\delta)$ to establish the desired result. □

*Proof of Theorem 3.2.* Let $Z_{t,n} = (X_t, t\lambda_n)$ and define the maps $h_\xi(Z_{t,n}) = 1\{X_t < \xi\} \times \cos(t\lambda_n)$ and $h_\xi^*(Z_{t,n}) = 1\{X_t < \xi\} \times \sin(t\lambda_n)$. The empirical process evaluated at some function $h$ is denoted by $\nu_n h := n^{-1/2} \sum_{t=1}^n (h(Z_{t,n}) - \mathbb{E}h(Z_{t,n}))$.

The finite Fourier transform at nonzero natural frequencies is invariant to centering. Hence, from (2.3), we have the decomposition

$$Q_{n,\tau}(\lambda_n) = \frac{1}{2\pi}\left|\nu_n(h_{\hat{\xi}_n(\tau)} - h_{\xi_0(\tau)}) + i\,\nu_n(h_{\hat{\xi}_n(\tau)}^* - h_{\xi_0(\tau)}^*) - n^{-1/2}\sum_{t=1}^n V_t(\tau)e^{-it\lambda_n}\right|^2 \quad \text{(A.1)}$$

For the proof of the theorem, I proceed in three steps: I show that (i) the first term and (ii) the second term inside the modulus in the display are small in probability and that (iii) the remainder of (A.1) has the desired asymptotic distribution jointly for frequencies $\lambda_n + 2\pi j/n$ with $|j| \leq k$.

*Step* (i): Define a norm $\rho(h_\xi) = \sup_{t,n \in \mathbb{N}} \|h_\xi(Z_{t,n})\|$. Take a grid of points $\xi_0(\tau) - \delta =: \xi_0 < \xi_1 < \cdots < \xi_N := \xi_0(\tau) + \delta$ and let $b_k(Z_{t,n}) := \big(h_{\xi_k}(Z_{t,n}) - h_{\xi_{k-1}}(Z_{t,n})\big)/\cos(t\lambda_n)$. Given a $\xi \in \mathcal{X}_\tau(\delta)$, we can then find an index $k$ such that $|h_\xi - h_{\xi_{k-1}}| \leq b_k$. In addition, we have

$$\rho(b_k) = \|1\{X_0 < \xi_k\} - 1\{X_0 < \xi_{k-1}\}\| \leq \sqrt{F_X(\xi_k) - F_X(\xi_{k-1})},$$

which is bounded above by a constant multiple of $\sqrt{\xi_k - \xi_{k-1}}$ due to Lipschitz continuity.

Hence, if we choose the grid such that $\rho(b_k) \leq \epsilon$ for all $k = 1, \ldots, N$, the parametric class $\mathcal{H} := \{h_\xi : \xi \in \mathcal{X}_\tau(\delta)\}$ has bracketing numbers (see Andrews and Pollard, 1994; van der Vaart, 1998, pp. 270-271) with respect to $\rho$ of order $N(\epsilon, \mathcal{H}) = O(\epsilon^{-2})$ as $\epsilon \to 0$.

By the same calculations as in the preceding display, there is some $M > 0$ such that all $\xi, \xi' \in \mathcal{X}_\tau(\delta)$ satisfy $\rho(h_\xi - h_{\xi'}) \leq M|\xi - \xi'|^{1/2}$ and therefore $\rho(h_{\hat{\xi}_n(\tau)} - h_{\xi_0(\tau)}) \to_p 0$ in view of Lemma A.1 below. For $\epsilon, \eta > 0$, the limit superior of $\mathbb{P}\big(|\nu_n(h_{\hat{\xi}_n(\tau)} - h_{\xi_0(\tau)})| \geq \epsilon\big)$ is then at most

$$\limsup_{n \to \infty} \mathbb{P}\Big(\big|\nu_n\big(h_{\hat{\xi}_n(\tau)} - h_{\xi_0(\tau)}\big)\big| \geq \epsilon, \rho\big(h_{\hat{\xi}_n(\tau)} - h_{\xi_0(\tau)}\big) \leq \eta\Big)$$

$$\leq \limsup_{n \to \infty} \mathbb{P}^*\Bigg(\sup_{\xi \in \mathcal{X}_\tau(\delta) : \rho(h_\xi - h_{\xi_0(\tau)}) \leq \eta} \big|\nu_n\big(h_\xi - h_{\xi_0(\tau)}\big)\big| \geq \epsilon\Bigg) \qquad (A.2)$$

The Markov inequality and Lemma A.3 below imply that the term on the right can be made as small as desired by choosing $\eta$ small enough. This is also true for the frequencies $\lambda_n + 2\pi j/n$ with $|j| \leq k$.

*Step* (ii): Replace cosines with sines in the proofs of Lemmas A.2 and A.3 (with the same bounding functions $b_k$ as above) to reach the same conclusion for $\nu_n(h^*_{\xi_0(\tau)} - h^*_{\hat{\xi}_n(\tau)})$.

*Step* (iii): In view of (i), (ii), and continuity of the modulus, I only have to show that the remainder of (A.1) converges jointly at each $\lambda_n + 2\pi j/n$, $|j| \leq k$, in distribution to independent exponential variables with mean $g_\tau(\lambda)$. For this I apply Corollary 2.1 of Shao and Wu (2007). Because $V_t(\tau)$ is a bounded mean-zero variable, the only condition that has to be checked is $\sum_{t=0}^{\infty} \|\mathbb{E}(V_t(\tau) \mid \mathcal{F}_0) - \mathbb{E}(V_t(\tau) \mid \mathcal{F}_{-1})\| < \infty$. By the conditional Jensen inequality, the law of iterated expectations, and Assumption C, this summability condition is satisfied because

$$\big\|\mathbb{E}\big(V_t(\tau) \mid \mathcal{F}_0\big) - \mathbb{E}\big(V_t(\tau) \mid \mathcal{F}_{-1}\big)\big\| = \big\|\mathbb{E}\big(1\{X_t < \xi_0(\tau)\} - 1\{X_t^* < \xi_0(\tau)\} \mid \mathcal{F}_0\big)\big\|$$
$$\leq \|1\{X_t < \xi_0(\tau)\} - 1\{X_t^* < \xi_0(\tau)\}\|.$$

Assumption A implies Assumption C, and so the joint convergence asserted in Theorem 3.2 follows. $\square$

**Lemma A.1.** *Suppose Assumptions B and C hold; then* $\sqrt{n}(\hat{\xi}_n(\tau) - \xi_0(\tau)) = O_p(1)$.

*Proof of Lemma A.1.* Arguing as in the proof of Theorem 1 of Wu (2005), use the conditional Jensen inequality and the law of iterated expectations to deduce that

$$\|\mathbb{E}(1\{X_t < \xi\} \mid \mathcal{F}_0) - \mathbb{E}(1\{X_t^* < \xi\} \mid \mathcal{F}_{-1}, \varepsilon_0^*)\|$$
$$= \|\mathbb{E}(1\{X_t < \xi\} - 1\{X_t^* < \xi\} \mid \mathcal{F}_0, \varepsilon_0^*)\|$$
$$\leq \|1\{X_t < \xi\} - 1\{X_t^* < \xi\}\|.$$

Taking suprema over $\mathcal{X}_\tau(\delta)$ shows that Assumption C implies condition (7) of Wu (2007) and his Theorem 1 then yields the desired result. $\square$

For Lemmas A.2 and A.3, I mimic the proofs of Andrews and Pollard's (1994) Theorem 2.2 and Lemma 3.1; their arguments do not apply directly since Andrews and Pollard work with strongly mixing arrays.

**Lemma A.2.** *Let $\phi(h_\xi - h_{\xi'}) := \rho(h_\xi - h_{\xi'})^{2/(2+\gamma)}$ for some $\gamma > 0$ and suppose that Assumption A holds. Then, for all $n \in \mathbb{N}$, all $\xi, \xi' \in \mathcal{X}_\tau(\delta)$, and every even integer $Q \geq 2$,*

$$\mathbb{E}|\nu_n(h_\xi - h_{\xi'})|^Q \leq n^{-Q/2} C\big((\phi(h_\xi - h_{\xi'})^2 n) + \cdots + (\phi(h_\xi - h_{\xi'})^2 n)^{Q/2}\big),$$

*where $C$ depends only on $Q$, $\gamma$, and $\sigma$. The inequality remains valid when $h_\xi - h_{\xi'}$ is replaced by $b_k$ for any given $k \geq 1$.*

*Proof of Lemma A.2.* It suffices to show the inequality given in the lemma after dividing both sides by $4^Q$ to ensure that the absolute value of

$$H_t := \big(h_\xi(Z_{t,n}) - h_{\xi'}(Z_{t,n}) - (\mathbb{E}h_\xi(Z_{t,n}) - \mathbb{E}h_{\xi'}(Z_{t,n}))\big)/4$$

is bounded by 1. The $4^{-Q}$ on the right hand can be absorbed into $C$. Define $H'_t$ in the same way as $H_t$ but replace $X_t$ with $X'_t$. Here I suppress the dependence of $H_t$ and $H'_t$ on $n$, $\xi$, and $\xi'$ because they are irrelevant in the following. Also note that $\mathbb{E}H_t = \mathbb{E}H'_t = 0$ for all $t, n \in \mathbb{N}$ and all $\xi, \xi' \in \mathcal{X}_\tau(\delta)$ because $X_t$ and $X'_t$ are identically distributed.

For fixed $k \geq 2$, $d \geq 1$, and $1 \leq m < k$, consider integers $t_1 \leq \cdots \leq t_m \leq t_{m+1} \leq \cdots \leq t_k$ so that $t_{m+1} - t_m = d$ and define $a_k(\lambda_n) = |2^{-k} \prod_{i=1}^k \cos(t_i \lambda_n)|$. Since $U_t := H_t/\cos(t\lambda_n)$ and $U'_t := H'_t/\cos(t\lambda_n)$ are stationary, repeatedly add and subtract to see that

$$\begin{aligned}
&\big|\mathbb{E}H_{t_1} H_{t_2} \cdots H_{t_k} - \mathbb{E}H_{t_1} H_{t_2} \cdots H_{t_m} \mathbb{E}H_{t_{m+1}} \cdots H_{t_k}\big| \\
&= a_k(\lambda_n)\big|\mathbb{E}U_{t_1-t_m} U_{t_2-t_m} \cdots U_{t_k-t_m} - \mathbb{E}U_{t_1-t_m} U_{t_2-t_m} \cdots U_0 \mathbb{E}U_d \cdots U_{t_k-t_m}\big| \\
&\leq a_k(\lambda_n)\big|\mathbb{E}U_{t_1-t_m} \cdots U_0 (U_d - U'_d) U_{t_{m+2}-t_m} \cdots U_{t_k-t_m}\big| \\
&\quad + \sum_{i=2}^{k-m-1} a_k(\lambda_n)\big|\mathbb{E}U_{t_1-t_m} \cdots U_0 U'_d \cdots (U_{t_{m+i}-t_m} - U'_{t_{m+i}-t_m}) \cdots U_{t_k-t_m}\big| \quad (\text{A.3}) \\
&\quad + a_k(\lambda_n)\big|\mathbb{E}U_{t_1-t_m} \cdots U_0 U'_d \cdots U'_{t_k-t_m} - \mathbb{E}U_{t_1-t_m} \cdots U_0 \mathbb{E}U_d \cdots U_{t_k-t_m}\big|
\end{aligned}$$

In particular, the last term on the right-hand side is zero because $U_{t_1-t_m} \cdots U_0$ and $U'_d \cdots U'_{t_k-t_m}$ are independent and $U'_d \cdots U'_{t_k-t_m}$ and $U_d \cdots U_{t_k-t_m}$ are identically distributed.

By Assumption A, $\|U_d - U'_d\|_s \leq \|1\{X_d < \xi\} - 1\{X'_d < \xi\}\|_s + \|1\{X_d < \xi'\} - 1\{X'_d < \xi'\}\|_s \leq 2\sup_{\xi \in \mathcal{X}_\tau(\delta)} \|1\{X_d < \xi\} - 1\{X'_d < \xi\}\|_s \leq C'\sigma^d$ for some $C' > 0$ and $s \geq 1$. Here the choice of $s$ does not matter because Assumption A still applies when $\|\cdot\|$ is replaced by $\|\cdot\|_s$ for any $s > 0$; see Lemma 2 of Wu and Min (2005). Hölder's inequality then bounds the first term on the right-hand side of the preceding display by

$$\|H_{t_1} \cdots H_{t_m}\|_p \|H_{t_{m+2}} \cdots H_{t_k}\|_q C' \sigma^d, \tag{A.4}$$

where the reciprocals of $p$, $q$, and $s$ sum to 1. Proceeding as in Andrews and Pollard (1994),

another application of the Hölder inequality yields

$$\|H_{t_1}\cdots H_{t_m}\|_p \leq \left(\prod_{i=1}^m \mathbb{E}|H_{t_i}|^{mp}\right)^{1/(mp)} \leq \phi(h_\xi - h_{\xi'})^{(2+\gamma)/p}$$

whenever $mp \geq 2$ and similarly $\|H_{t_{m+2}}\cdots H_{t_k}\|_q \leq \phi(h_\xi - h_{\xi'})^{(2+\gamma)/q}$ whenever $(k - m - 1)q \geq 2$. Suppose for now that $k \geq 3$. If $k > m + 1$, take $s = (\gamma + Q)/\gamma$ and $mp = (k - m - 1)q = (k - 1)/(1 - 1/s)$. Decrease the resulting exponent of $\phi(h_\xi - h_{\xi'})$ from $Q(2 + \gamma)/(Q + \gamma)$ to 2 to see that (A.4) is bounded by $C'\sigma^d\phi(h_\xi - h_{\xi'})^2$. If $k \geq 2$ and $k = m + 1$, the factor $\|H_{t_{m+2}}\cdots H_{t_k}\|_q$ is not present in (A.4), but we can still choose $s = (\gamma + Q)/\gamma$ and $mp = (k - 1)/(1 - 1/s)$ to obtain the same bound. Since the same argument also applies to each of the other summands in (A.3), we can find a constant $M > 0$ such that

$$\left|\mathbb{E}H_{t_1}H_{t_2}\cdots H_{t_k}\right| \leq \left|\mathbb{E}H_{t_1}H_{t_2}\cdots H_{t_m}\mathbb{E}H_{t_{m+1}}\cdots H_{t_k}\right| + M\sigma^d\phi(h_\xi - h_{\xi'})^2.$$

Here $M$ in fact depends on $k$, but this does not disturb any of the subsequent steps.

Now replace (A.2) in Andrews and Pollard (1994) by the inequality in the preceding display. In particular, replace their $8\alpha(d)^{1/s}$ with $M\sigma^d$ and their $\tau^2$ with $\phi(h_\xi - h_{\xi'})^2$. The rest of their arguments now go through without changes.

The inequality for $b_k$ follows by letting $\lambda_n \equiv 0$; this is not a contradiction to the assumptions of Theorem 3.2 because this proof is valid for any sequence $(\lambda_n)_{n\in\mathbb{N}}$. □

**Lemma A.3.** *Suppose the assumptions of Theorem 3.2 hold. For every $\epsilon > 0$ and every even integer $Q \geq 4$, there is an $\eta > 0$ such that*

$$\limsup_{n\to\infty} \mathbb{E}^*\left(\sup_{\xi,\xi'\in\mathcal{X}_\tau(\delta):\rho(h_\xi-h_{\xi'})\leq\eta} |\nu_n(h_\xi - h_{\xi'})|\right)^Q \leq \epsilon.$$

*Proof of Lemma A.3.* I follow Andrews and Pollard's (1994) proof of their Theorem 2.1. It requires three steps: (i) Their "Proof of inequality (3.2)," (ii) their "Proof of inequality (3.3)," and (iii) their "Comparison of pairs" argument. Replace their $i$ with $k$ and their $\tau(h_i)$ with $\phi(b_k)$; then apply Lemma A.2 above instead of Andrews and Pollard's (1994) Lemma 3.1 in the derivation of their inequality (3.5) to deduce

$$\left\|\max_{1\leq k\leq N}|\nu_n b_k|\right\|_Q \leq C'N^{1/Q}\max\left\{n^{-1/2}, \max_{1\leq k\leq N}\phi(b_k)\right\}$$

and use this inequality in (i) instead of their inequality (3.5). In (i) Andrews and Pollard also require the finiteness of the bracketing integral $\int_0^1 x^{-\gamma/(2+\gamma)}N(x,\mathcal{H})^{1/Q}\,dx$, which follows immediately by choosing $\gamma = Q - 4$. Another application of Lemma A.2 establishes the required analogue of Andrews and Pollard's inequality (3.5) used in (ii). The same inequality can also be applied in (iii). The other arguments remain valid without changes. □

*Proof of Theorem 3.6.* Denote by $\tilde{r}_{n,\tau}(j) = n^{-1}\sum_{t=|j|+1}^{n} V_t(\tau)V_{t-|j|}(\tau)$ and

$$g_{n,\tau}(\lambda) := \frac{1}{2\pi}\sum_{|j|<n} w(j/B_n)\tilde{r}_{n,\tau}(j)e^{-ij\lambda}$$

the infeasible sample auto-covariance and smoothed quantile spectrum, respectively, based on the unknown quantile $\xi_0(\tau)$. The triangle inequality and $|V_t(\cdot,\cdot)| < 1$ yield

$$2\pi \sup_{\lambda \in (-\pi,\pi]} |\hat{g}_{n,\tau}(\lambda) - g_{n,\tau}(\lambda)|$$

$$\leq \frac{1}{n}\sum_{|j|<n}|w(j/B_n)|\sum_{t=|j|+1}^{n}|\hat{V}_t(\tau)\hat{V}_{t-|j|}(\tau) - V_t(\tau)V_{t-|j|}(\tau)|$$

$$\leq \frac{1}{n}\sum_{|j|<n}|w(j/B_n)|\sum_{t=|j|+1}^{n}\left(|\hat{V}_t(\tau) - V_t(\tau)| + |\hat{V}_{t-|j|}(\tau) - V_{t-|j|}(\tau)|\right)$$

$$\leq \frac{1}{n}\sum_{|j|<n}|w(j/B_n)|\sum_{t=|j|+1}^{n}\left(\mathbf{1}_{\{|X_t - \xi_0(\tau)| < |\hat{\xi}_n(\tau) - \xi_0(\tau)|\}}\right.$$

$$\left. + \mathbf{1}_{\{|X_{t-|j|} - \xi_0(\tau)| < |\hat{\xi}_n(\tau) - \xi_0(\tau)|\}}\right).$$

Consider the first indicator function on the right-hand side of the preceding display and recall that $\sqrt{n}(\hat{\xi}_n(\tau) - \xi_0(\tau))$ is uniformly tight by Lemma A.1. For any given $\epsilon > 0$ and $\eta > 0$, the Markov inequality implies for large enough $M > 0$

$$\limsup_{n\to\infty}\mathbb{P}\left(\frac{1}{n}\sum_{|j|<n}|w(j/B_n)|\sum_{t=|j|+1}^{n}\mathbf{1}_{\{|X_t - \xi_0(\tau)| < |\hat{\xi}_n(\tau) - \xi_0(\tau)|\}} \geq \eta\right)$$

$$\leq \limsup_{n\to\infty}\mathbb{P}\left(\frac{1}{n}\sum_{|j|<n}|w(j/B_n)|\sum_{t=|j|+1}^{n}\mathbf{1}_{\{|X_t - \xi_0(\tau)| < Mn^{-1/2}\}} \geq \eta\right)$$

$$+ \sup_{n\in\mathbb{N}}\mathbb{P}(|\hat{\xi}_n(\tau) - \xi_0(\tau)| \geq Mn^{-1/2})$$

$$\leq \limsup_{n\to\infty}\frac{1}{\eta}\sum_{|j|<n}|w(j/B_n)|\mathbb{P}(|X_0 - \xi_0(\tau)| < Mn^{-1/2}) + \epsilon.$$

By Lemma 1 of Jansson (2002), the limit superior of $B_n^{-1}\sum_{|j|<n}|w(j/B_n)|$ is finite, and in view of the assumed Lipschitz continuity, the first term on the right-hand side of the preceding display then vanishes because $B_n n^{-1/2} \to 0$. The same argument applies to the second indicator function above due to stationarity. Together this yields $\sup_\lambda |\hat{g}_{n,\tau}(\lambda) - g_{n,\tau}(\lambda)| \to_p 0$.

To show $g_{n,\tau}(\lambda) \to_p g_\tau(\lambda)$ uniformly in $\lambda$, I use Liu and Wu's (2010) Theorem 1, which applies whenever the windows $w \in \mathcal{W}$ satisfy their Condition 1. The only two conditions that need to be established are the absolute integrability of $w$, which is immediate from $\int_{-\infty}^{\infty} |w(x)| \, dx \leq 2 \int_0^{\infty} \bar{w}(x) \, dx < \infty$, and

$$\limsup_{n \to \infty} B_n^{-1} \sum_{j \in \mathbb{Z}} w(j/B_n)^2 < \infty.$$

Although Liu and Wu (2010) provide a specific value for the limit in the preceding display, its boundedness is in fact all that is needed for the proof of their Theorem 1. To this end, take $M \geq \sup_{x \in \mathbb{R}} |w(x)|$ such that for $j \geq 1$

$$w(j/B_n)^2 \leq M|w(j/B_n)| \leq M\bar{w}(j/B_n) \leq MB_n \int_{(j-1)/B_n}^{j/B_n} \bar{w}(x) \, dx$$

by monotonicity, and therefore symmetry implies

$$B_n^{-1} \sum_{j \in \mathbb{Z}} w(j/B_n)^2 \leq B_n^{-1} + 2M \sum_{j=1}^{\infty} \int_{(j-1)/B_n}^{j/B_n} \bar{w}(x) \, dx = B_n^{-1} + 2M \int_0^{\infty} \bar{w}(x) \, dx,$$

which is finite by assumption. This is also true for its limit superior as $n \to \infty$. The triangle inequality completes the proof. $\square$

*Proof of Theorem 3.7.* Note that $\sqrt{m_n/B_n}(\tilde{g}_{m_n,\tau}(\lambda) - \mathbb{E}g_{m_n,\tau}(\lambda))$ can be written as

$$\sqrt{\frac{m_n}{B_n}} \frac{1}{2\pi} \sum_{|j|<m_n} w(j/B_n) e^{-ij\lambda} \frac{1}{m_n} \sum_{t=1+m_n}^{n-|j|} \left( V_t(\tau, \tilde{\xi}_n(\tau)) V_{t+|j|}(\tau, \tilde{\xi}_n(\tau)) - \mathbb{E}V_t(\tau)V_{t+|j|}(\tau) \right).$$

Let $\tilde{V}_t(\tau, \xi) = F_{\tilde{X}}(\xi_0(\tau)) - 1\{\tilde{X}_t < \xi\}$ and $\tilde{V}_t(\tau) = F_{\tilde{X}}(\xi_0(\tau)) - 1\{\tilde{X}_t < \xi_0(\tau)\}$. Decompose $V_t(\tau, \tilde{\xi}_n(\tau))V_{t+|j|}(\tau, \tilde{\xi}_n(\tau)) - \mathbb{E}V_t(\tau)V_{t+|j|}(\tau)$ into

$$\tilde{V}_t(\tau, \tilde{\xi}_n(\tau))\tilde{V}_{t+|j|}(\tau, \tilde{\xi}_n(\tau)) - \mathbb{E}\tilde{V}_t(\tau, \tilde{\xi}_n(\tau))\tilde{V}_{t+|j|}(\tau, \tilde{\xi}_n(\tau)) \quad (A.5)$$

$$- \left( \tilde{V}_t(\tau)\tilde{V}_{t+|j|}(\tau) - \mathbb{E}\tilde{V}_t(\tau)\tilde{V}_{t+|j|}(\tau) \right) \quad (A.6)$$

$$+ V_t(\tau, \tilde{\xi}_n(\tau))V_{t+|j|}(\tau, \tilde{\xi}_n(\tau)) - \tilde{V}_t(\tau, \tilde{\xi}_n(\tau))\tilde{V}_{t+|j|}(\tau, \tilde{\xi}_n(\tau)) \quad (A.7)$$

$$+ \mathbb{E}\tilde{V}_t(\tau, \tilde{\xi}_n(\tau))\tilde{V}_{t+|j|}(\tau, \tilde{\xi}_n(\tau)) - \mathbb{E}\tilde{V}_t(\tau)\tilde{V}_{t+|j|}(\tau) \quad (A.8)$$

$$+ \mathbb{E}\tilde{V}_t(\tau)\tilde{V}_{t+|j|}(\tau) - \mathbb{E}V_t(\tau)V_{t+|j|}(\tau) \quad (A.9)$$

$$+ \tilde{V}_t(\tau)\tilde{V}_{t+|j|}(\tau) - \mathbb{E}\tilde{V}_t(\tau)\tilde{V}_{t+|j|}(\tau), \quad (A.10)$$

where the expectations are evaluated for a given $\tilde{\xi}_n(\tau)$. The idea is to show that the contributions of (A.5)-(A.9) to the spectral density estimate are negligible in large samples. Then I apply central limit theory of Shao and Wu (2007) to the smoothed quantile periodogram constructed from (A.10).

(A.5) and (A.6): Let $f_\xi(\tilde{X}_t, \tilde{X}_{t+|j|}) = \tilde{V}_t(\tau, \xi)\tilde{V}_{t+|j|}(\tau, \xi)$ and $p_n = p_n(j) = \lfloor (n-|j|)/(|j|+l_n) \rfloor$. Up to a remainder term that will be negligible in the following, we can write

$$\sum_{t=1+m_n}^{n-|j|} \left(f_\xi(\tilde{X}_t, \tilde{X}_{t+|j|}) - \mathbb{E}f_\xi(\tilde{X}_0, \tilde{X}_{|j|})\right)$$

$$= \sqrt{p_n} \sum_{k=1+m_n}^{|j|+l_n} p_n^{-1/2} \sum_{t=1}^{p_n} \left(f_\xi(\tilde{X}_{k+(t-1)(|j|+l_n)}, \tilde{X}_{k+(t-1)(|j|+l_n)+|j|}) - \mathbb{E}f_\xi(\tilde{X}_0, \tilde{X}_{|j|})\right).$$

Take $\tilde{\nu}_{p_n,k}$ to be the empirical process constructed from the $p_n$ observations contained in the inner sums on the right. The preceding display can then be rewritten as

$$\sqrt{p_n} \sum_{k=1+m_n}^{|j|+l_n} \tilde{\nu}_{p_n,k} f_\xi.$$

Notice that each of the empirical processes is a sum of iid random variables. Consider the class $\mathcal{F}_\eta = \{f_\xi - f_{\xi_0(\tau)} : |\xi - \xi_0(\tau)| \leq \eta\}$ for some $\eta > 0$ with envelope

$$F_\eta(\tilde{X}_t, \tilde{X}_{t+|j|}) = 1\{|\tilde{X}_t - \xi_0(\tau)| \leq \eta\} + 1\{|\tilde{X}_{t+|j|} - \xi_0(\tau)| \leq \eta\}.$$

The subgraph of a real-valued function $f$ is defined as $\{(x,t) : f(x) > t\}$. Write the subgraph of $(f_\xi - f_{\xi_0(\tau)})(x_1, x_2)$ as

$(\{x_1 < \xi\} \cap \{x_2 < \xi\} \cap \{x_1 < \xi_0(\tau)\} \cap \{x_2 < \xi_0(\tau)\} \cap \{(f_\xi - f_{\xi_0(\tau)})(x_1, x_2) > t\})$
$\cup (\{x_1 \geq \xi\} \cap \{x_2 < \xi\} \cap \{x_1 < \xi_0(\tau)\} \cap \{x_2 < \xi_0(\tau)\} \cap \{(f_\xi - f_{\xi_0(\tau)})(x_1, x_2) > t\})$
$\cup (\{x_1 < \xi\} \cap \{x_2 \geq \xi\} \cap \{x_1 < \xi_0(\tau)\} \cap \{x_2 < \xi_0(\tau)\} \cap \{(f_\xi - f_{\xi_0(\tau)})(x_1, x_2) > t\})$
$\cup \cdots \cup$
$\cup (\{x_1 \geq \xi\} \cap \{x_2 \geq \xi\} \cap \{x_1 \geq \xi_0(\tau)\} \cap \{x_2 \geq \xi_0(\tau)\} \cap \{(f_\xi - f_{\xi_0(\tau)})(x_1, x_2) > t\}).$

Note that within each member of the union in the display, the function $(f_\xi - f_{\xi_0(\tau)})(x_1, x_2)$ is constant and therefore a Vapnik–Červonenkis (VC) class of sets. More generally, the collection of sets $\{z \in \mathbb{R}^2 : z^\top \lambda \leq 0\}$ indexed by $\lambda \in \mathbb{R}^2$ is a VC class (see van der Vaart and Wellner, 1996, Problem 2.6.14); the same holds for the collection $\{z \in \mathbb{R}^2 : z^\top \lambda > 0\}$ by Lemma 2.6.17(i) of van der Vaart and Wellner. Each remaining set in the preceding display indexed by $\xi \in \mathbb{R}$ and $t \in \mathbb{R}$ is contained in one of these two VC classes. Subclasses of VC classes are VC classes themselves. Conclude from van der Vaart and Wellner's Lemma 2.6.17(ii) and (iii) that the subgraph of $(f_\xi - f_{\xi_0(\tau)})(x_1, x_2)$ is a VC class. Therefore $\mathcal{F}_\eta$ is a VC subgraph class. As such, it satisfies Pollard's (1982) uniform entropy condition; for a modern statement of this condition see van der Vaart and Wellner (1996, p. 239, with the condition being $J(1, \mathcal{F}_\eta) < \infty$ in their notation).

Use van der Vaart and Wellner's (1996) Theorem 2.14.1 to deduce that there exists some absolute constant $C > 0$ so that

$$\left\| \sup_{\xi:|\xi-\xi_0(\tau)|\leq \eta} |\tilde{\nu}_{p_n,k}(f_\xi - f_{\xi_0(\tau)})| \right\|_2 \leq C\|F_\eta\|_2.$$

Lipschitz continuity implies that the right-hand side can be bounded further by a constant multiple of $\sqrt{\eta}$, uniformly in $j$. Pick $\eta = Mn^{-1/2}$ for some large $M > 0$. Then

$$\mathbb{P}\left(\left|(m_nB_n)^{-1/2}\frac{1}{2\pi}\sum_{|j|<m_n} w(j/B_n)e^{-ij\lambda}\sqrt{p_n}\sum_{k=1+m_n}^{|j|+l_n} \tilde{\nu}_{p_n,k}(f_{\tilde{\xi}_n(\tau)} - f_{\xi_0(\tau)})\right| \geq \epsilon\right)$$

$$\leq \mathbb{P}\left((m_nB_n)^{-1/2}\frac{1}{2\pi}\sum_{|j|<m_n} |w(j/B_n)|\sqrt{p_n}\sum_{k=1+m_n}^{|j|+l_n} \sup_{\xi:|\xi-\xi_0(\tau)|\leq Mn^{-1/2}} |\tilde{\nu}_{p_n,k}(f_\xi - f_{\xi_0(\tau)})| \geq \epsilon\right)$$

$$+ \sup_{n\in\mathbb{N}} \mathbb{P}\big(|\tilde{\xi}_n(\tau) - \xi_0(\tau)| \geq Mn^{-1/2}\big).$$

By the Markov inequality and the tail bound derived above, the square of the first term is bounded by a constant multiple of

$$(m_nB_n)^{-1/2} \sum_{|j|<m_n} |w(j/B_n)|\sqrt{p_n} \sum_{k=1+m_n}^{|j|+l_n} \left\| \sup_{\xi:|\xi-\xi_0(\tau)|\leq Mn^{-1/2}} |\tilde{\nu}_{p_n,k}(f_\xi - f_{\xi_0(\tau)})| \right\|_2,$$

which is $O(B_n/n^{1/4}) + O(\sqrt{B_nl_n}/n^{1/4}) = o(1)$. Since $M$ can be chosen to be as large as desired without affecting the other bounds, conclude that the contribution of (A.5) and (A.6) is negligible.

(A.7): The absolute value of (A.7) is at most

$$\left|V_t(\tau, \tilde{\xi}_n(\tau)) - \tilde{V}_t(\tau, \tilde{\xi}_n(\tau))\right| + \left|V_{t+|j|}(\tau, \tilde{\xi}_n(\tau)) - \tilde{V}_{t+|j|}(\tau, \tilde{\xi}_n(\tau))\right|$$
$$\leq 2|\tau - F_{\tilde{X}}(\xi_0(\tau))| + |\mathbf{1}\{X_t < \tilde{\xi}_n(\tau)\} - \mathbf{1}\{\tilde{X}_t < \tilde{\xi}_n(\tau)\}|$$
$$+ |\mathbf{1}\{X_{t+|j|} < \tilde{\xi}_n(\tau)\} - \mathbf{1}\{\tilde{X}_{t+|j|} < \tilde{\xi}_n(\tau)\}| \tag{A.11}$$

By the triangle inequality and the GMC property, the first term on the right-hand side is bounded by $2\mathbb{E}|\mathbf{1}\{X_t < \xi_0(\tau)\} - \mathbf{1}\{\tilde{X}_t < \xi_0(\tau)\}| = O(\varrho^{l_n}) \leq O(n^{-1})$. Because $t > \lfloor n/2 \rfloor$, the expectation of the second term can be bounded by

$$\mathbb{P}(|\tilde{X}_t - \tilde{\xi}_n(\tau)| \leq |X_t - \tilde{X}_t|) \leq \mathbb{P}(|\tilde{X}_t - \tilde{\xi}_n(\tau)| \leq \varrho^{l_n/2}) + \mathbb{P}(|X_t - \tilde{X}_t| > \varrho^{l_n/2})$$
$$\leq \mathbb{E}\big(F_{\tilde{X}}(\tilde{\xi}_n(\tau) + \varrho^{l_n/2}) - F_{\tilde{X}}(\tilde{\xi}_n(\tau) - \varrho^{l_n/2})\big) + O(\varrho^{l_n/2}),$$

where the second inequality follows from the law of iterated expectations, independence of $\tilde{X}_t$ and $\tilde{\xi}_n(\tau)$, and the GMC property. By Lipschitz continuity, the right-hand side of the preceding display is $O(\varrho^{l_n/2}) \leq O(n^{-1})$. The same arguments apply to (A.11). These results hold uniformly in $j$, so

$$\sqrt{\frac{m_n}{B_n}} \mathbb{E}\left|\frac{1}{2\pi} \sum_{|j|<m_n} w(j/B_n)e^{-ij\lambda}\frac{1}{m_n} \sum_{t=1+m_n}^{n-|j|} \left(V_t(\tau,\tilde{\xi}_n(\tau))V_{t+|j|}(\tau,\tilde{\xi}_n(\tau))\right.\right.$$
$$\left.\left. - \tilde{V}_t(\tau,\tilde{\xi}_n(\tau))\tilde{V}_{t+|j|}(\tau,\tilde{\xi}_n(\tau))\right)\right|$$
$$\leq \sqrt{\frac{m_n}{B_n}} \sum_{|j|<m_n} |w(j/B_n)|O(n^{-1}) = O(\sqrt{B_n/n}) = o(1).$$

Negligibility now follows from the Markov inequality.

(A.8): For $|j| \geq l_n$, $\mathbb{E}\tilde{V}_t(\tau)\tilde{V}_{t+|j|}(\tau) = \mathbb{E}\tilde{V}_t(\tau)\mathbb{E}\tilde{V}_{t+|j|}(\tau) = 0$ and (A.8) equals

$$\mathbb{E}\tilde{V}_t(\tau,\tilde{\xi}_n(\tau))\mathbb{E}\tilde{V}_{t+|j|}(\tau,\tilde{\xi}_n(\tau)) = \left(F_{\tilde{X}}(\xi_0(\tau)) - F_{\tilde{X}}(\tilde{\xi}_n(\tau))\right)^2 = O_p(n^{-1})$$

uniformly in $j$ by stationarity and Lipschitz continuity. For $|j| < l_n$, the absolute value of (A.8) is does not exceed two times

$$\mathbb{E}|1(\tilde{X}_t < \tilde{\xi}_n(\tau)) - 1(\tilde{X}_t < \xi_0(\tau))| \leq \mathbb{P}(|\tilde{X}_t - \xi_0(\tau)| \leq |\tilde{\xi}_n(\tau) - \xi_0(\tau)|) = O_p(n^{-1/2})$$

uniformly in $j$. (Recall that the expectation is with respect to the randomness in $\tilde{X}_t$.) Conclude that

$$\sqrt{\frac{m_n}{B_n}}\left|\frac{1}{2\pi}\sum_{|j|<m_n} w(j/B_n)e^{-ij\lambda}\frac{1}{m_n}\sum_{t=1+m_n}^{n-|j|}\left(\mathbb{E}\tilde{V}_t(\tau,\tilde{\xi}_n(\tau))\tilde{V}_{t+|j|}(\tau,\tilde{\xi}_n(\tau)) - \mathbb{E}\tilde{V}_t(\tau)\tilde{V}_{t+|j|}(\tau)\right)\right|$$
$$\leq \sqrt{\frac{m_n}{B_n}}\sum_{|j|<l_n}|w(j/B_n)|O_p(n^{-1/2}) + \sqrt{\frac{m_n}{B_n}}\sum_{l_n\leq|j|<m_n}|w(j/B_n)|O_p(n^{-1})$$
$$= O_p(l_n B^{-1/2}) + O_p(\sqrt{B_n/n}) = o_p(1).$$

(A.9): The absolute value of (A.9) is, uniformly in $j$, bounded by two times

$$\mathbb{E}|\tilde{V}_t(\tau) - V_t(\tau)| \leq |\tau - F_{\tilde{X}}(\xi_0(\tau))| + \mathbb{E}|1\{X_t < \xi_0\} - 1\{\tilde{X}_t < \xi_0\}|$$
$$\leq 2\mathbb{E}|1\{X_t < \xi_0(\tau)\} - 1\{\tilde{X}_t < \xi_0(\tau)\}| = O(\varrho^{l_n}) \leq O(n^{-1}),$$

so the contribution to the spectral estimate is negligible for the same reasons as in (A.7).

(A.10): The proof of Theorem 3.1 of Shao and Wu (2007) immediately implies that

$$\sqrt{\frac{m_n}{B_n}}\frac{1}{2\pi}\sum_{|j|<m_n} w(j/B_n)e^{-ij\lambda}\frac{1}{m_n}\sum_{t=1+m_n}^{n-|j|}\left(\tilde{V}_t(\tau)\tilde{V}_{t+|j|}(\tau) - \mathbb{E}\tilde{V}_t(\tau)\tilde{V}_{t+|j|}(\tau)\right)$$

has an asymptotic normal distribution with the variance given in the statement of the theorem. The sequence $(l_n)_{n\in\mathbb{N}}$ used here differs slightly from the one given in Shao and Wu's proof of their Theorem 3.1, but inspection of their proof reveals that this does not disturb any of their arguments. The bias term $\mathbb{E}g_{m_n,\tau}(\lambda) - g_\tau(\lambda)$ is $O(B_n^{-3})$ by a routine argument and the assumed GMC property. $\square$

*Proof of Theorem 4.1.* (i) The process $S_{n,\tau}(\lambda)$ can be decomposed into

$$\sqrt{n}\sum_{j=1}^{n-1}(\hat{r}_{n,\tau}(j) - \tilde{r}_{n,\tau}(j))\psi_j(\lambda) + \sqrt{n}\sum_{j=1}^{n-1}\tilde{r}_{n,\tau}(j)\psi_j(\lambda).$$

The second term side converges weakly in $L_2(\Pi)$ to $S_\tau(\lambda)$ by the proof of Theorem 2.1 of Shao (2011a). The Continuous Mapping Theorem then yields $CM_{n,\tau} \rightsquigarrow \|S_\tau\|_\Pi^2$ as long as the $L_2(\Pi)$-norm of the first term of the display is eventually small in probability.

To this end, define $\varphi_j = \|\psi_j\|_\Pi$ and note that $\langle \psi_j, \psi_k \rangle = 0$ for $j \neq k$. Use this orthogonality to write

$$\left\|\sqrt{n}\sum_{j=1}^{n-1}(\hat{r}_{n,\tau}(j) - \tilde{r}_{n,\tau}(j))\psi_j\right\|_\Pi^2 = n\sum_{j=1}^{n-1}(\hat{r}_{n,\tau}(j) - \tilde{r}_{n,\tau}(j))^2 \varphi_j^2.$$

Let $r_{n,\tau}(j,\xi) = (\tau - F_X(\xi))^2(n-j)/n$ for $j > 0$ and $\xi \in \mathbb{R}$. Under the null hypothesis, we have $r_\tau(j) = 0 = r_{n,\tau}(j,\xi_0(\tau))$ for all $j > 0$ and, by Lipschitz continuity, there exists $M > 0$ such that $\check{r}_{n,\tau}(j) := r_{n,\tau}(j,\hat{\xi}_n(\tau)) \leq |F_X(\hat{\xi}_n(\tau)) - \tau|^2 \leq M|\hat{\xi}_n(\tau) - \xi_0(\tau)|^2$ for all $j > 0$. In view of these properties, fix some $K \leq n$ and apply the Loève $c_r$ inequality to bound the quantity in the preceding display by

$$2n \sum_{j=1}^{K-1}\left(\hat{r}_{n,\tau}(j) - \check{r}_{n,\tau}(j) - (\tilde{r}_{n,\tau}(j) - r_\tau(j))\right)^2 \varphi_j^2 \tag{A.12}$$

$$+ 2n \sum_{j=1}^{K-1}\check{r}_{n,\tau}(j)^2 \varphi_j^2 + n\sum_{j=K}^{n-1}(\hat{r}_{n,\tau}(j) - \tilde{r}_{n,\tau}(j))^2 \varphi_j^2. \tag{A.13}$$

Fix $\epsilon, \epsilon' > 0$ and let $e_{\xi,j}(X_t, X_{t+j}) := V_t(\tau,\xi)V_{t+j}(\tau,\xi)$. For given $j > 0$, take $\rho(e_{\xi,j} - e_{\xi',j})$ as the distance of $\xi$ and $\xi'$ on $\mathcal{X}_\tau(\delta)$, where $\rho$ is as in the proof of Theorem 3.2. The distance also depends on $\tau$, but this is irrelevant in the following. Note that $\rho(e_{\xi,j} - e_{\xi',j}) \leq 2\|1\{X_0 < \xi\} - 1\{X_0 < \xi'\}\| \leq M'|\xi - \xi'|$ uniformly in $j > 0$ for some $M' > 0$ by stationarity. Hence, for any $\eta_j > 0$, $1 \leq j < K$, we have

$$\mathbb{P}\bigcup_{j=1}^{K-1}\left\{\rho\left(e_{\hat{\xi}_n(\tau),j} - e_{\xi_0(\tau),j}\right) > \eta_j\right\} \leq \mathbb{P}\left(|\hat{\xi}_n(\tau) - \xi_0(\tau)|^{1/2} > \min_{1\leq j<K}\eta_j/M'\right) = o(1).$$

Under the null hypothesis, we can write $\sqrt{n}(\tilde{r}_{n,\tau}(j) - r_\tau(j)) = \sqrt{(n-j)/n}\,\nu_{n-j}\,e_{\xi_0(\tau)}$ and, by (4.3), $\sqrt{n}(\hat{r}_{n,\tau}(j) - \check{r}_{n,\tau}(j)) = \sqrt{(n-j)/n}\,\nu_{n-j}\,e_{\hat{\xi}_n(\tau)}$ as long as $\hat{\xi}_n(\tau) \in \mathcal{X}_\tau(\delta)$, where I use the notation from the proof of Theorem 3.2. In view of the preceding display, (A.12) then satisfies

$$\limsup_{n\to\infty} \mathbb{P}\left(n\sum_{j=1}^{K-1}\left(\hat{r}_{n,\tau}(j) - \check{r}_{n,\tau}(j) - (\tilde{r}_{n,\tau}(j) - r_\tau(j))\right)^2 \varphi_j^2 \geq \epsilon/2\right)$$

$$\leq \limsup_{n\to\infty} \mathbb{P}^*\left(\sum_{j=1}^{K-1}\left(\sup_{\xi \in \mathcal{X}_\tau(\delta):\rho(e_{\xi,j}-e_{\xi_0(\tau),j})\leq \eta_j} |\nu_{n-j}(e_{\xi,j} - e_{\xi_0(\tau),j})|\right)^2 \varphi_j^2 \geq \epsilon/2\right)$$

$$\leq \frac{2}{\epsilon}\sum_{j=1}^{K-1} \limsup_{n\to\infty} \mathbb{E}^*\left(\sup_{\xi \in \mathcal{X}_\tau(\delta):\rho(e_{\xi,j}-e_{\xi_0(\tau),j})\leq \eta_j} |\nu_{n-j}(e_{\xi,j} - e_{\xi_0(\tau),j})|\right)^2 \varphi_j^2$$

$$\leq \frac{2}{\epsilon}\sum_{j=1}^{K-1}\left(\limsup_{n\to\infty} \mathbb{E}^*\left(\sup_{\xi \in \mathcal{X}_\tau(\delta):\rho(e_{\xi,j}-e_{\xi_0(\tau),j})\leq \eta_j} |\nu_{n-j}(e_{\xi,j} - e_{\xi_0(\tau),j})|\right)^Q\right)^{2/Q} \varphi_j^2$$

$$\leq \frac{2\epsilon'}{\epsilon}\sum_{j=1}^{K-1} \varphi_j^2 \leq \frac{2\epsilon'}{\epsilon}\sum_{j>0} \varphi_j^2 = \left(\frac{\pi}{6\epsilon}\right)\epsilon',$$

where the first inequality is the Markov inequality, the second follows from Lemma A.4 below, the third is Jensen's, and the equality uses $\varphi_j^2 = 1/(2\pi j^2)$ for $j > 0$.

Now consider (A.13). The first term can be bounded by

$$2nM^2|\hat{\xi}_n(\tau) - \xi_0(\tau)|^4 \sum_{j>0} \varphi_j^2 = O_p(n^{-1}) = o_p(1).$$

By (4.3) and arguments as in the proof of Theorem 3.6, for a large enough $M'$ the probability that the second term of (A.13) exceeds $\epsilon$ is at most

$$\mathbb{P}\left(n\sum_{j=K}^{n-1}\left(\hat{r}_{n,\tau}(j) - \tilde{r}_{n,\tau}(j)\right)^2 \varphi_j^2 \geq \epsilon, |\hat{\xi}_n(\tau) - \xi_0(\tau)| \leq M'n^{-1/2}\right) + \epsilon'$$

$$\leq \mathbb{P}\left(n^{-1}\sum_{j=K}^{n-1}\left(\sum_{t=j+1}^{n} 1_{\{|X_t-\xi_0(\tau)|<|\hat{\xi}_n(\tau)-\xi_0(\tau)|\}}\right.\right.$$

$$\left.\left. + 1_{\{|X_{t-j}-\xi_0(\tau)|<|\hat{\xi}_n(\tau)-\xi_0(\tau)|\}}\right)^2 \varphi_j^2 \geq \epsilon, |\hat{\xi}_n(\tau) - \xi_0(\tau)| \leq M'n^{-1/2}\right) + \epsilon'$$

$$\leq \mathbb{P}\left(n^{-1}\sum_{j=K}^{n-1}\left(\sum_{t=j+1}^{n} 1_{\{|X_t-\xi_0(\tau)|\leq M'n^{-1/2}\}}\right)^2 \varphi_j^2\right.$$

$$\left. + n^{-1}\sum_{j=K}^{n-1}\left(\sum_{t=j+1}^{n} 1_{\{|X_{t-j}-\xi_0(\tau)|\leq M'n^{-1/2}\}}\right)^2 \varphi_j^2 \geq \epsilon/2\right) + \epsilon'$$

$$\leq \frac{4}{\epsilon}\mathbb{E}\left(n^{-1/2}\sum_{t=1}^{n}1_{\{|X_t-\xi_0(\tau)|\leq M'n^{-1/2}\}}\right)^2\sum_{j\geq K}\varphi_j^2 + \epsilon'$$

$$\leq \frac{4}{\epsilon}\left(\mathbb{P}(|X_0-\xi_0(\tau)|\leq M'n^{-1/2}) + n\mathbb{P}(|X_0-\xi_0(\tau)|\leq M'n^{-1/2})^2\right)\sum_{j\geq K}\varphi_j^2 + \epsilon'$$

$$= O(1)\sum_{j\geq K}\varphi_j^2 + \epsilon',$$

which can be made smaller than $2\epsilon'$ by choosing $K$ large enough. This does not affect any of the other bounds. Since $\epsilon'$ was arbitrary, we have $CM_{n,\tau} = \tilde{CM}_{n,\tau} + o_p(1)$ and $CM_{n,\tau} \rightsquigarrow \|S_\tau\|_\Pi^2$, which proves the first result.

(ii) Fix some $K \leq n$ and decompose the statistic into

$$CM_n(\tau)/n = \sum_{j=1}^{K-1}(\hat{r}_{n,\tau}(j) - \tilde{r}_{n,\tau}(j))^2\varphi_j^2 + 2\sum_{j=1}^{K-1}(\hat{r}_{n,\tau}(j) - \tilde{r}_{n,\tau}(j))\tilde{r}_{n,\tau}(j)\varphi_j^2$$
$$+ \sum_{j=1}^{K-1}\tilde{r}_{n,\tau}(j)^2\varphi_j^2 + \sum_{j=K}^{n-1}\hat{r}_{n,\tau}(j)^2\varphi_j^2.$$

The first and second terms on the right-hand side of the displayed equation converge to zero in probability as $n \to \infty$ because $\hat{r}_{n,\tau}(j) - \tilde{r}_{n,\tau}(j) = o_p(1)$ for each $j$ under the assumptions of the theorem. The third term converges in probability to $\sum_{j=1}^{K-1}r_\tau(j)^2\varphi_j^2$ as $n \to \infty$ by Wu's (2005) Theorem 2(i). The absolute value of the last term is bounded by $\sum_{j\geq K}\varphi_j^2$, where I have used the fact that $|\hat{r}_{n,\tau}| \leq 1$. Hence, let $K \to \infty$ to conclude $CM_n(\tau)/n \to_p \sum_{j>0}r_\tau(j)^2\varphi_j^2 > 0$. The desired result now follows from a routine argument. $\square$

**Lemma A.4.** *Suppose the assumptions of Theorem 4.1 hold. For every $j > 0$, every $Q$, and every $\epsilon > 0$, there is an $\eta > 0$ such that*

$$\limsup_{n\to\infty}\mathbb{E}^*\left(\sup_{\xi,\xi'\in\mathcal{X}_\tau(\delta):\rho(e_{\xi,j}-e_{\xi',j})\leq\eta}|\nu_{n-j}(e_{\xi,j} - e_{\xi',j})|\right)^Q \leq \epsilon.$$

*Proof of Lemma A.4.* Take a grid of points $\xi_0(\tau) - \delta = \xi_0 < \xi_1 < \cdots < \xi_N = \xi_0(\tau) + \delta$ and let $b_{k,j}(X_t, X_{t+j}) := 1\{X_t < \xi_k\} - 1\{X_t < \xi_{k-1}\} + 1\{X_{t+j} < \xi_k\} - 1\{X_{t+j} < \xi_{k-1}\}$. Given a $\xi \in \mathcal{X}_\tau(\delta)$, we can then find an index $k$ such that $|e_{\xi,j} - e_{\xi_{k-1},j}| \leq b_{k,j}$. Further,

$$\rho(b_{k,j}) \leq 2\|1\{X_0 < \xi_k\} - 1\{X_0 < \xi_{k-1}\}\| \leq 2\sqrt{F_X(\xi_k) - F_X(\xi_{k-1})},$$

which is proportional to $\sqrt{\xi_k - \xi_{k-1}}$ due to Lipschitz continuity. If $\rho(b_{k,j}) \leq \epsilon$ for all $k = 1, \ldots, N$ then, as above, for each $j$ the parametric class $\mathcal{E}_j := \{e_{\xi,j} : \xi \in \mathcal{X}_\tau(\delta)\}$ has bracketing numbers with respect to $\rho$ of order $N(\epsilon, \mathcal{E}_j) = O(\epsilon^{-2})$ as $\epsilon \to 0$. Hence, bracketing integrals of the class $\mathcal{H}$ above and the classes $\mathcal{E}_j$ have the same behavior. The proof of Lemma A.3 therefore also applies to this lemma as long as the reference to Lemma A.2 is replaced by Lemma A.5 below. $\square$

**Lemma A.5.** *Fix some $\gamma > 0$ and suppose that Assumption A holds. For all $n \in \mathbb{N}$, all $j < n$, all $\xi, \xi' \in \mathcal{X}_\tau(\delta)$, and every even integer $Q \geq 2$ we have*

$$\mathbb{E}|\nu_{n-j}(e_{\xi,j} - e_{\xi',j})|^Q \leq (n-j)^{-Q/2} C\big((\phi(e_{\xi,j} - e_{\xi',j})^2 (n-j)) \\ + \cdots + (\phi(e_{\xi,j} - e_{\xi',j})^2 (n-j))^{Q/2}\big),$$

*where $C$ depends only on $j$, $Q$, $\gamma$, and $\sigma$. The inequality remains valid when $e_{\xi,j} - e_{\xi',j}$ is replaced by $b_{k,j}$ for any given $k \geq 1$.*

*Proof of Lemma A.5.* As in the proof of Lemma A.2, it suffices to show the inequality given in the Lemma after dividing both sides by $4^Q$ to ensure that the absolute value of

$$E_{t,t+j} := \big(e_{\xi,j}(X_t, X_{t+j}) - e_{\xi',j}(X_t, X_{t+j}) - (\mathbb{E}e_{\xi,j}(X_0, X_j) - \mathbb{E}e_{\xi',j}(X_0, X_j))\big)/4$$

is bounded by 1. Define $E'_{t,t+j}$ in the same way as $E_{t,t+j}$ but replace $X_t$ with $X'_t$ and $X_{t+j}$ with $X'_{t+j}$. For fixed $k \geq 2$, $d \geq 1$, and $1 \leq m < k$, consider integers $t_1 \leq \cdots \leq t_m \leq t_{m+1} \leq \cdots \leq t_k$ so that $t_{m+1} - t_m = d$. Repeatedly add and subtract to see that

$$\big|\mathbb{E}E_{t_1,t_1+j} \cdots E_{t_k,t_k+j} - \mathbb{E}E_{t_1,t_1+j} \cdots E_{t_m,t_m+j}\mathbb{E}E_{t_{m+1},t_{m+1}+j} \cdots E_{t_k,t_k+j}\big|$$
$$= \big|\mathbb{E}E_{t_1-t_m-j,t_1-t_m} \cdots E_{t_k-t_m-j,t_k-t_m}$$
$$\qquad - \mathbb{E}E_{t_1-t_m-j,t_1-t_m} \cdots E_{-j,0}\mathbb{E}E_{l-j,l} \cdots E_{t_k-t_m-j,t_k-t_m}\big|$$
$$\leq \big|\mathbb{E}E_{t_1-t_m-j,t_1-t_m} \cdots E_{-j,0}(E_{d-j,d} - E'_{d-j,d})$$
$$\qquad \times E_{t_{m+2}-t_m-j,t_{m+2}-t_m} \cdots E_{t_k-t_m-j,t_k-t_m}\big|$$
$$+ \sum_{i=2}^{k-m-1} \big|\mathbb{E}E_{t_1-t_m-j,t_1-t_m} \cdots E_{-j,0}E'_{d-j,d} \times \cdots$$
$$\qquad \times (E_{t_{m+i}-t_m-j,t_{m+i}-t_m} - E'_{t_{m+i}-t_m-j,t_{m+i}-t_m}) \cdots E_{t_k-t_m-j,t_k-t_m}\big|$$
$$+ \big|\mathbb{E}E_{t_1-t_m-j,t_1-t_m} \cdots E_{-j,0}E'_{d-j,d} \cdots E'_{t_k-t_m-j,t_k-t_m}$$
$$\qquad - \mathbb{E}E_{t_1-t_m-j,t_1-t_m} \cdots E_{-j,0}\mathbb{E}E_{l-j,l} \cdots E_{t_k-t_m-j,t_k-t_m}\big|,$$

where the last term on the right-hand side can again been seen to be zero.

Since $j$ is fixed, it is possible to write $\|E_{d-j,d} - E'_{d-j,d}\|_s \leq \|1\{X_d < \xi\} - 1\{X'_d < \xi\}\|_s + \|1\{X_{d-j} < \xi'\} - 1\{X'_{d-j} < \xi'\}\|_s \leq C'\sigma^d(1 + \sigma^{-j})$ for some $C' > 0$, where the cases where $d \leq j$ were absorbed into $C'$. The same can then be done for $(1 + \sigma^{-j})$. Hence, proceed exactly as above to find a constant $M > 0$ such that

$$\big|\mathbb{E}E_{t_1,t_1+j}E_{t_2,t_2+j} \cdots E_{t_k,t_k+j}\big|$$
$$\leq \big|\mathbb{E}E_{t_1,t_1+j}E_{t_2,t_2+j} \cdots E_{t_m,t_m+j}\mathbb{E}E_{t_{m+1},t_{m+1}+j} \cdots E_{t_k,t_k+j}\big| + M\sigma^d \phi(e_{\xi,j} - e_{\xi',j})^2.$$

The rest of the arguments in the proof of Lemma A.2 now go through without changes. The proof for the bounding functions $b_k$ is almost identical and therefore omitted. □

*Proof of Corollary 4.5.* (i) Theorem 1 of Lifshits (1982) guarantees that $\|S_\tau\|_\Pi^2$ has a con-

tinuous distribution function, and therefore $c_{n,\tau}(1-\alpha) \to c_{\infty,\tau}(1-\alpha)$ by Lemma 21.2 of van der Vaart (1998), where $c_{\infty,\tau}$ is the quantile function of $\|S_\tau\|_\Pi^2$. Hence, $CM_{n,\tau} - c_{n,\tau}(1-\alpha) \rightsquigarrow \|S_\tau\|_\Pi^2 - c_{\infty,\tau}(1-\alpha)$ and, in particular, $\|S_\tau\|_\Pi^2 - c_\infty(1-\alpha)$ also has a continuous distribution function. This in turn implies

$$\bigl|\mathbb{P}\bigl(CM_{n,\tau} > c_{n,\tau}(1-\alpha)\bigr) - \alpha\bigr| = \bigl|\mathbb{P}\bigl(CM_{n,\tau} \leq c_{n,\tau}(1-\alpha)\bigr) - \mathbb{P}\bigl(\|S_\tau\|_\Pi^2 \leq c_{\infty,\tau}(1-\alpha)\bigr)\bigr| \to 0.$$

(ii) Let $CM_{\infty,\tau} := \sum_{j>0} r_\tau(j)^2 \varphi_j^2$ and pick an $\epsilon > 0$ such that $CM_{\infty,\tau} - \epsilon > 0$. By Theorem 4.1(ii) and the properties of quantile functions,

$$\mathbb{P}\bigl(CM_{n,\tau} \leq c_{n,\tau}(1-\alpha)\bigr) \leq 1\bigl(c_{n,\tau}(1-\alpha) > n(CM_{\infty,\tau} - \epsilon)\bigr)$$
$$+ \mathbb{P}(|CM_{n,\tau}/n - CM_{\infty,\tau}| \geq \epsilon)$$
$$= 1\bigl(1 - \alpha > \mathbb{P}(CM'_{n,\tau}/n \leq CM_{\infty,\tau} - \epsilon)\bigr) + o(1).$$

It therefore suffices to show that $CM'_{n,\tau}/n \to_p 0$, which follows from an application of Birkhoff's Ergodic Theorem to the first term on the right-hand side of

$$n^{-1}CM'_{n,\tau} \leq \sum_{j=1}^{K-1}\left(n^{-1}\sum_{t=1+j}^{n}(\tau - J_t)(\tau - J_{t-j})\right)^2 \varphi_j^2 + \sum_{j\geq K} \varphi_j^2$$

and then letting $K \to \infty$. $\square$

*Proof of Theorem 4.7.* (i) Recall that $r_\tau(j) = 0$ for all $j > 0$ under the null hypothesis and let $\tilde{r}_{n,\tau}^*(j) := n^{-1}\sum_{t=j+1}^{n} V_t(\tau)V_{t-j}(\tau)\omega_t$. We can write $CM_{n,\tau}^* = \|S_{n,\tau}^*\|_\Pi^2$, where $S_{n,\tau}^*$ is

$$n^{-1/2}\sum_{j=1}^{n-1}\left(\sum_{t=j+1}^{n}\bigl(\hat{V}_t(\tau)\hat{V}_{t-j}(\tau) - V_t(\tau)V_{t-j}(\tau)\bigr)\omega_t\right)\psi_j(\lambda) \qquad (A.14)$$

$$- n^{-1/2}\sum_{j=1}^{n-1}\hat{r}_{n,\tau}(j)\psi_j(\lambda)\left(\sum_{t=j+1}^{n}\omega_t\right) + \sqrt{n}\sum_{j=1}^{n-1}\tilde{r}_{n,\tau}^*(j)\psi_j(\lambda). \qquad (A.15)$$

As a preliminary step, I show that the $L_2(\Pi)$-norms the first two terms have a $\hat{\mathbb{P}}$-probability limit of zero with high $\mathbb{P}$-probability; the $L_2(\Pi)$-norm of the third term converges $\hat{\mathbb{P}}$-weakly in $\mathbb{P}$-probability to $\|S_\tau\|_\Pi^2$ by Shao's (2011a) Theorem 3.1. I then use these results below to prove that the bootstrap test has asymptotic size $\alpha$.

The $\hat{\mathbb{P}}$-expectation of the square of the $L_2(\Pi)$-norm of (A.14) can be written as

$$\hat{\mathbb{E}}\left\|n^{-1/2}\sum_{j=1}^{n-1}\left(\sum_{t=j+1}^{n}\bigl(\hat{V}_t(\tau)\hat{V}_{t-j}(\tau) - V_t(\tau)V_{t-j}(\tau)\bigr)\omega_t\right)\psi_j\right\|_\Pi^2$$
$$= n^{-1}\sum_{j=1}^{n-1}\varphi_j^2 \sum_{s=1}^{L_n}\left(\sum_{t\in\mathcal{B}_s \cap [j+1,n]}\bigl(\hat{V}_t(\tau)\hat{V}_{t-j}(\tau) - V_t(\tau)V_{t-j}(\tau)\bigr)\right)^2.$$

Fix $\epsilon, \epsilon' > 0$ and pick a large enough $M > 0$ such that $\sup_{n \in \mathbb{N}} \mathbb{P}(|\hat{\xi}_n(\tau) - \xi_0(\tau)| > Mn^{-1/2}) < \epsilon'$. As in the proof of Theorem 4.1(i), the probability that the term on the right is larger than $\epsilon$ is at most $\epsilon'$ plus

$$\frac{2}{\epsilon} \sum_{j=1}^{n-1} \varphi_j^2 n^{-1} \sum_{s=1}^{L_n} \mathbb{E}\left(\sum_{t \in \mathcal{B}_s \cap [j+1,n]} 1_{\{|X_t - \xi_0(\tau)| \leq Mn^{-1/2}\}}\right)^2$$

$$\leq \frac{2}{\epsilon} \sum_{j=1}^{n-1} \varphi_j^2 b_n^{-1} \left(b_n \mathbb{P}(|X_0 - \xi_0(\tau)| \leq Mn^{-1/2}) + b_n^2 \mathbb{P}(|X_0 - \xi_0(\tau)| \leq Mn^{-1/2})^2\right)$$

$$\leq \frac{2}{\epsilon}\left(O(n^{-1/2}) + O(b_n/n)\right) \sum_{j>0} \varphi_j^2,$$

which can be made arbitrarily small by first letting $n \to \infty$ and then $M \to \infty$.

Now consider the $\hat{\mathbb{P}}$-expectation of the square of the $L_2(\Pi)$-norm of (A.15), which can be written as

$$\hat{\mathbb{E}}\left\|n^{-1/2} \sum_{j=1}^{n-1} \hat{r}_{n,\tau}(j)\left(\sum_{t=j+1}^{n} \omega_t\right)\psi_j\right\|_{\Pi}^2 = n^{-1} \sum_{j=1}^{n-1} \hat{r}_{n,\tau}(j)^2 \varphi_j^2 \sum_{s=1}^{L_n} \hat{\mathbb{E}}\left(\sum_{t \in \mathcal{B}_s \cap [j+1,n]} \omega_t\right)^2$$

$$\leq b_n \sum_{j=1}^{n-1} \hat{r}_{n,\tau}(j)^2 \varphi_j^2$$

$$\leq 2b_n \sum_{j=1}^{n-1} \left(\hat{r}_{n,\tau}(j) - \tilde{r}_{n,\tau}(j)\right)^2 \varphi_j^2 + 2b_n \sum_{j=1}^{n-1} \tilde{r}_{n,\tau}(j)^2 \varphi_j^2$$

by the Loève $c_r$ inequality. The first term on the right-hand side of the display converges to zero in probability by arguments similar to those given in the proof of Theorem 4.1(i) provided that $b_n/n \to 0$. The second term is $O_p(b_n/n)$ by Corollary 2.1 of Shao (2011a). It follows that $\sup_x |\hat{\mathbb{P}}(CM_{n,\tau}^* \leq x) - \mathbb{P}(\|S_\tau\|_\Pi^2 \leq x)| = o_p(1)$.

Theorem 1 of Lifshits (1982) and Lemma 21.2 of van der Vaart (1998) then give $c_{n,\tau}^*(1-\alpha) \to_p c_\infty(1-\alpha)$. Thus, $CM_{n,\tau} - c_{n,\tau}^*(1-\alpha) \rightsquigarrow \|S_\tau\|_\Pi^2 - c_{\infty,\tau}(1-\alpha)$, which yields

$$\left|\mathbb{P}(CM_{n,\tau} > c_{n,\tau}^*(1-\alpha)) - \alpha\right|$$
$$= \left|\mathbb{P}(CM_{n,\tau} \leq c_{n,\tau}^*(1-\alpha)) - \mathbb{P}(\|S_\tau\|_\Pi^2 \leq c_{\infty,\tau}(1-\alpha))\right| \to 0.$$

(ii) Recall that $CM_{\infty,\tau} = \sum_{j>0} r_\tau(j)^2 \varphi_j^2$. Pick an $\epsilon > 0$ such that $CM_{\infty,\tau} > \epsilon$ and, as in the proof of Corollary 4.5(ii), the properties of quantile functions and Theorem 4.1(ii) imply

$$\mathbb{P}(CM_{n,\tau} \leq c_{n,\tau}^*(1-\alpha)) \leq \mathbb{P}(c_{n,\tau}^*(1-\alpha) > n(CM_{\infty,\tau} - \epsilon))$$
$$+ \mathbb{P}(|CM_{n,\tau}/n - CM_{\infty,\tau}| \geq \epsilon)$$
$$= \mathbb{P}(1-\alpha > \hat{\mathbb{P}}(CM_{n,\tau}^*/n \leq CM_{\infty,\tau} - \epsilon)) + o(1).$$

Hence it suffices to show that $\hat{\mathbb{E}} CM^*_{n,\tau}/n \to_p 0$, which is seen from

$$n^{-1}\hat{\mathbb{E}} CM^*_{n,\tau} = n^{-2}\sum_{j=1}^{n-1}\varphi_j^2 \sum_{s=1}^{L_n}\left(\sum_{t\in\mathcal{B}_s\cap[j+1,n]}\left(\hat{V}_t(\tau)\hat{V}_{t+j}(\tau) - \hat{r}_{n,\tau}(j)\right)\right)^2 \leq 4\frac{b_n}{n}\sum_{j>0}\varphi_j^2$$

almost surely and $b_n/n \to 0$. $\square$

## References

Andrews, D. W. K. (1991). Heteroskedasticity and autocorrelation consistent covariance matrix estimation. *Econometrica 59*, 817–858.

Andrews, D. W. K. and D. Pollard (1994). An introduction to functional central limit theorems for dependent stochastic processes. *International Statistical Review 62*, 119–132.

Basraka, B., R. A. Davis, and T. Mikosch (2002). Regular variation of garch processes. *Stochastic Processes and their Applications 99*, 95–115.

Bollerslev, T. (1986). Generalized autoregressive conditional heteroskedasticity. *Journal of Econometrics 31*, 307–327.

Bougerol, P. and N. Picard (1992). Strict stationarity of generalized autoregressive processes. *Annals of Probability 20*, 1714–1730.

Brillinger, D. R. (1975). *Time Series: Data Analysis and Theory*. Holt, Rinehart and Winston.

Brockwell, P. J. and R. A. Davis (1991). *Time Series: Theory and Methods* (2nd ed.). Springer, New York.

Chernozhukov, V., C. Hansen, and M. Jansson (2009). Finite sample inference for quantile regression models. *Journal of Econometrics 152*, 93–103.

Chung, J. and Y. Hong (2007). Model-free evaluation of directional predictability in foreign exchange markets. *Journal of Applied Econometrics 22*, 855–889.

Cogley, T. and T. J. Sargent (2002). Evolving post-World War II U.S. inflation dynamics. In *NBER Macroeconomics Annual 2001*, Volume 17, pp. 331–388.

Cont, R. (2001). Empirical properties of asset returns: Stylized facts and statistical issues. *Quantitative Finance 1*, 223–236.

Davidson, J. and R. de Jong (2000). Consistency of kernel estimators of heteroscedastic and autocorrelated covariance matrices. *Econometrica 68*, 407–424.

Davidson, J., A. Monticini, and D. Peel (2007). Implementing the wild bootstrap using a two-point distribution. *Economics Letters 96*, 309–315.

Dette, H., M. Hallin, T. Kley, and S. Volgushev (2011). Of copulas, quantiles, ranks and spectra: an $L_1$-approach to spectral analysis. Unpublished manuscript, `arXiv:1111.7205v1`.

Diebold, F. X., L. E. Ohanian, and J. Berkowitz (1998). Dynamic equilibrium economies: A framework for comparing models and data. *Review of Economic Studies 65*, 433–451.

Ding, Z., C. Granger, and R. F. Engle (1993). A long memory property of stock market returns and a new model. *Journal of Empirical Finance 1*, 83–106.

Durbin, J. (1967). Tests of serial independence based on the cumulated periodogram. *Bulletin of the International Statistical Institute 42*, 1040–1048.

Engle, R. F. (1982). Autoregressive conditional heteroscedasticity with estimates of the variance of United Kingdom inflation. *Econometrica 50*, 987–1007.

Giacomini, R., D. N. Politis, and H. White (2007). A warp-speed method for conducting Monte Carlo experiments involving bootstrap estimators. Unpublished manuscript, Department of Economics, University of California, Los Angeles.

Granger, C. W. J. (1966). The typical spectral shape of an economic variable. *Econometrica 34*, 150–161.

Hagemann, A. (2011). Robust spectral analysis. Unpublished manuscript, Department of Economics, University of Illinois, `arXiv:1111.1965v1`.

Heyde, C. C. (2002). On modes of long-range dependence. *Journal of Applied Probability 39*, 882–888.

Hong, Y. (1996). Consistent testing for serial correlation of unknown form. *Econometrica 64*, 837–864.

Hong, Y. (1999). Hypothesis testing in time series via the empirical characteristic function: A generalized spectral density approach. *Journal of the American Statistical Association 94*, 1201–1220.

Hsing, T. and W. B. Wu (2004). On weighted $U$-statistics for stationary processes. *Annals of Probability 32*, 1600–1631.

Huber, P. J. and E. M. Ronchetti (2009). *Robust Statistics* (Second ed.). John Wiley & Sons.

Jansson, M. (2002). Consistent covariance matrix estimation for linear processes. *Econometric Theory 18*, 1449–1459.

Katkovnik, V. (1998). Robust $M$-periodogram. *IEEE Transactions on Signal Processing 46*, 3104–3109.

Kleiner, B. and R. D. Martin (1979). Robust estimation of power spectra. *Journal of the Royal Statistical Society. Series B. 49*, 313–351.

Klüppelburg, C. and T. Mikosch (1994). Some limit theory for the self-normalised periodogram of stable processes. *Scandinavian Journal of Statistics 21*, 485–491.

Knight, K. (2006). Comment on quantile autoregression. *Journal of the American Statistical Association 101*, 994–996.

Koenker, R. W. (2005). *Quantile Regression*. Cambridge University Press, New York.

Koenker, R. W. and G. Bassett (1978). Regression quantiles. *Econometrica 46*, 33–50.

Koenker, R. W. and Z. Xiao (2006). Quantile autoregression. *Journal of the American Statistical Association 101*, 980–1006.

Lee, J. and S. Subba Rao (2011). The quantile spectral density and comparison based tests for nonlinear time series. Unpublished manuscript, Department of Statistics, Texas A&M University `arXiv:1112.2759v1`.

Li, T.-H. (2008). Laplace periodogram for time series analysis. *Journal of the American Statistical Association 103*, 757–768.

Li, T.-H. (2012). Quantile periodograms. *Journal of the American Statistical Association 107*, 765–776.

Lifshits, M. A. (1982). On the absolute continuity of distributions of functionals of random processes. *Theory of Probability and Its Applications 27*, 600–607.

Ling, S. and M. McAleer (2002). Stationarity and the existence of moments of a family of GARCH processes. *Journal of Econometrics 106*, 109–117.

Linton, O. and J.-Y. Whang (2007). The quantilogram: With an application to evaluating directional predictability. *Journal of Econometrics 141*, 250–282.

Liu, R. Y. (1988). Bootstrap procedures under some non-I.I.D. models. *Annals of Statistics 16*, 1696–1708.

Liu, W. and W. B. Wu (2010). Asymptotics of spectral density estimates. *Econometric Theory 26*, 1218–1248.

Longin, F. M. (1996). The asymptotic distribution of extreme stock market returns. *Journal of Business 69*, 383–408.

Loretan, M. and P. C. B. Phillips (1994). Testing the covariance stationarity of heavy-tailed time series: An overview of the theory with applications to several financial datasets. *Journal of Empirical Finance 1*, 211–248.

Mammen, E. (1992). *When Does Bootstrap Work? Asymptotic Results and Simulations*, Volume 77 of *Lecture Notes in Statistics*. Springer, New York.

Milhøj, A. (1981). A test of fit in time series models. *Biometrika 68*, 177–187.

Parzen, E. (1957). On consistent estimates of the spectrum of a stationary time series. *Annals of Mathematical Statistics 29*, 329–348.

Pollard, D. (1982). A central limit theorem for empirical processes. *Journal of the Australian Mathematical Mathematical Society (Series A) 33*, 235–248.

Priestley, M. B. (1962). Basic considerations in the estimation of spectra. *Technometrics 4*, 551–564.

Qu, Z. and D. Tkachenko (2011). Identification and frequency domain QML estimation of linearized DSGE models. Unpublished manuscript, Department of Economics, Boston University.

Robinson, P. M. (1991). Automatic frequency domain inference on semiparametric and nonparametric models. *Econometrica 59*, 1329–1363.

Rosenblatt, M. (1984). Asymptotic normality, strong mixing, and spectral density estimates. *Annals of Probability 12*, 1167–1180.

Rosenblatt, M. and J. W. Van Ness (1965). Estimation of the bispectrum. *Annals of Mathematical Statistics 36*, 1120–1136.

Sargent, T. J. (1987). *Macroeconomic Theory*. Academic Press, New York.

Shao, X. (2011a). A bootstrap-assisted spectral test of white noise under unknown dependence. *Journal of Econometrics 162*, 213–224.

Shao, X. (2011b). Testing for white noise under unknown dependence and its applications to goodness-of-fit for time series models. Forthcoming in *Econometric Theory*.

Shao, X. and W. B. Wu (2007). Asymptotic spectral theory for nonlinear time series. *Annals of Statistics 35*, 1773–1801.

Smith, R. J. (2005). Automatic positive semidefinite HAC covariance matrix and GMM estimation. *Econometric Theory 21*, 158–170.

van der Vaart, A. W. (1998). *Asymptotic Statistics*. Cambridge University Press.

van der Vaart, A. W. and J. A. Wellner (1996). *Weak Convergence and Empirical Processes: With Applications to Statistics*. Springer Series in Statistics. Springer, New York.

Walsh, J. E. (1960). Nonparametric tests for median by interpolation from sign tests. *Annals of the Institute of Statistical Mathematics 11*, 183–188.

Wu, W. B. (2005). Nonlinear system theory: Another look at dependence. *Proceedings of the National Academy of Sciences 102*, 14150–14154.

Wu, W. B. (2007). $M$-estimation of linear models with dependent errors. *Annals of Statistics 35*, 495–521.

Wu, W. B. and W. Min (2005). On linear processes with dependent observations. *Stochastic Processes and their Applications 115*, 939–958.